%% file: manuscript.tex
\documentclass[onefignum,onetabnum]{siamart190516}

\graphicspath{{figures/}{./}}
\input{ex_shared}

\input{abbrv}

\ifpdf
\hypersetup{
  pdftitle={A neural network approach for stochastic optimal control},
  pdfauthor={Xingjian Li, Deepanshu Verma and Lars Ruthotto}
}
\fi

\begin{document}

\maketitle

% REQUIRED
\begin{abstract}
  We present a neural network approach for approximating the value function of high-dimensional stochastic control problems. Our training process simultaneously updates our value function estimate and identifies the part of the state space likely to be visited by optimal trajectories.
   Our approach leverages insights from optimal control theory and the fundamental relation between semi-linear parabolic partial differential equations and forward-backward stochastic differential equations. 
   To focus the sampling on relevant states during neural network training, we use the stochastic Pontryagin maximum principle (PMP) to obtain the optimal controls for the current value function estimate.
  By design, our approach coincides with the method of characteristics for the non-viscous Hamilton-Jacobi-Bellman equation arising in deterministic control problems.
  Our training loss consists of a weighted sum of the objective functional of the control problem and penalty terms that enforce the HJB equations along the sampled trajectories.
  Importantly, training is unsupervised in that it does not require solutions of the control problem.
  
  Our numerical experiments highlight our scheme's ability to identify the relevant parts of the state space and produce meaningful value estimates.
  Using a two-dimensional model problem, we demonstrate the importance of the stochastic PMP to inform the sampling and compare to a finite element approach. With a nonlinear control affine quadcopter example, we illustrate that our approach can handle complicated dynamics.
  For a 100-dimensional benchmark problem, we demonstrate that our approach improves accuracy and time-to-solution and, via a modification, we show the wider applicability of our scheme.
\end{abstract}

% REQUIRED

% REQUIRED

\section{Introduction}
We consider stochastic optimal control (SOC) problems that require finding a policy to control randomly perturbed dynamical systems to optimize a given objective functional.
Problems of this type arise in many fields, including finance, biology, robotics, and many other engineering applications; see, for example, \cite{flemingsoner06,Stengel-1994} for references and extensive theoretical discussion on SOC problems.

Dynamic programming (DP) is a prominent framework for solving SOC problems.
At its core, DP seeks to find the value function, which assigns every state of the system the optimal cost-to-go.
For many problems, recovering the optimal control from the value function is straightforward.
One way to compute the value function is by solving the Hamilton-Jacobi-Bellman (HJB) equation \cite{Bellman58-DPSOC,yong_zhou-1999}, which is a semi-linear parabolic Partial Differential Equation (PDE).
The two main challenges in solving the HJB equations are the forward-backward structure and the high dimensionality.

The space dimension of the HJB equation equals the state dimension of the dynamical system to be controlled. 
Hence, many numerical schemes for solving PDEs cannot be applied to realistic problem instances due to the Curse-of-Dimensionality (CoD); for example, the computational costs of approaches that discretize the state space with a mesh typically grow exponentially in the dimensions of the system. 
Hence, the applicability of existing approaches that employ spatial discretizations (for example, \cite{Bertsekas95Tsitsiklis,Dong-Krylov-2007,Jakobsen03,Krylov05,Powell-book,Kushner-1990,num_methods_kushner_dupuis-2001,Powell-book,central_diff_HJB_wang_forsyth}) is limited to state dimensions $d \leq 3$. Mesh-free methods based on radial basis function methods (for example, \cite{RBF-FD}) can be effective for slightly larger $d$ but ultimately also suffer from CoD for $d = \mathcal{O}(100)$, as we consider here.

A first requirement for mitigating the CoD in SOC problems is to alleviate the need for spatial discretization.
This can be achieved by using a nonlinear version of the Feynman-Kac lemma and replacing the HJB equation with a system of Forward-Backward Stochastic Differential Equations (FBSDEs); see, for example, \cite{Cheridito_Touzi_2007,Pardoux_Tang_1999,yong_zhou-1999}.
Discretizing the FBSDE system using the classical Euler Maruyama scheme \cite{Kloeden-1992} yields mesh-free numerical schemes that have gained importance for SOC problems and non-linear second-order PDEs; see, for example, \cite{Bender,Exarchos-2019,Exarchos-2018,Panos-2020,Ma-Yong} and references therein.

A second requirement for mitigating the CoD is to parameterize the value function effectively in high dimensions.
Due to their universal approximation property and many advances in deep learning, the use of neural networks (NNs) as function approximators has recently gathered significant attention. 
The idea of combining the BSDE approach with NNs has been pioneered by the seminal works \cite{Han-2017,Han_BSDE} that use neural networks to solve a wider class of high-dimensional semi-linear parabolic PDEs and include a stochastic optimal control example as a special case.
Around the same time, ~\cite{raissi2018forward} achieves similarly promising results for neural-network-based solution of the HJB equation.
These works developed learning algorithms for approximating the value function or its gradient using the FBSDE system. 
Their success has sparked several extensions that improve the NN architectures and consider other types of PDEs; see, for example, \cite{Pham-2020,Peng-2020,pereira20a}.

Even after combining FBSDEs and deep learning, solving the HJB equation globally remains affected by the curse of dimensionality as it would require sampling the entire state space.
Fortunately, the optimal state trajectories of many SOC problems do not cover the entire space. This motivates us to develop semi-global approaches; that is, we seek to reliably estimate the value function in those parts of the state space that are likely to be visited when following an optimal policy.
Clearly, these tasks are interrelated: Finding the relevant parts of the state space depends on the value function, and, vice versa, the value function needs to be trained using samples from the state space.

The key idea of our approach is to use the stochastic Pontryagin maximum principle (PMP) \cite{pontryagin62} to link the sampling and the value function approximation; more precisely, we define the forward system in the FBSDE approach in terms of the value function. 
The PMP provides a set of necessary optimality conditions, known as the Hamiltonian system, and a closed loop feedback form that allows recovering the optimal control from the value function; see, for example, \cite{yong_zhou-1999}.  Our formulation is suitable for problems where the underlying Hamiltonian can be efficiently computed, such as those involving affine controls and convex Lagrangians with differentiable value functions. As we discuss in more detail later, the choice of the forward SDE is crucial as it is used to sample the parts of the state space along which the value function estimate is improved using a loss that includes the backward SDE.
It is worth noting that the FBSDE framework allows choosing the forward system almost arbitrarily.
Thus, our choice is theoretically on par with the standard Brownian motion used for exploring the state space in~\cite{Han_BSDE, raissi2018forward}. However, as our numerical experiments suggest, our choice of the forward dynamics enables us to solve a wider class of SOC problems; see~\cref{sec:num_exp}.

A theoretical advantage of our proposed scheme over existing approaches is its consistency with deterministic optimal control problems.
In the absence of uncertainty in the dynamics, our approach reduces to the method of characteristics for the HJB equation. 
Hence, our approach can also be seen as an extension of the neural network approaches for deterministic optimal control proposed, for example, in \cite{onken2020neural,Karl-Walter-2021}.

The rest of the paper is organized as follows:
In \cref{s:prelim}, we introduce the parts of SOC theory used to obtain our proposed forward SDE.
In \cref{s:NN-approach}, we provide the NN approach to learn the value function using the FBSDE reformulation from \cref{s:prelim}. Using an intuitive 2D example, we show the importance of sampling. 
For this low-dimensional problem, we use a finite element approach to compare and validate our method. 
To illustrate the potential of our method, we consider the 100-dimensional benchmark problem also used in~\cite{Han_BSDE,raissi2018forward}, and propose modifications that lead to more challenging sampling problems. We also test our method on a 12-dimensional problem with complex nonlinear dynamics.
Finally, we conclude the paper and discuss future directions.

%
%%%%%%%%%%%%%%%%%%%%%%%%%%%%%%%%%%%%%%%%%%%%%%%%%%%%%%%%
\section{Stochastic Optimal Control Background}
\label{s:prelim}
%%%%%%%%%%%%%%%%%%%%%%%%%%%%%%%%%%%%%%%%%%%%%%%%%%%%%%%%
In this section, we describe the stochastic optimal control (SOC) problems considered in this work and review the key results from SOC theory that motivate our approach.
Our discussion follows \cite{yong_zhou-1999} and we refer to this textbook for a more comprehensive background and for more general results. 
We first introduce the SOC problem and then review its underlying theory; to be specific, we review the stochastic Pontryagin Maximum Principle (PMP), the Hamilton-Jacobi-Bellman (HJB) equation, and its reformulation into a system of forward-backward stochastic differential equations (FBSDEs) obtained from a nonlinear version of the Feynman-Kac formula.

%%%%%%%%%%%%%%%%%%%%%%%%%%%%%%%%%%%%%%%%%%%
\subsection{Stochastic Optimal Control Problem}
%%%%%%%%%%%%%%%%%%%%%%%%%%%%%%%%%%%%%%%%%%%
Let $(\Omega, \Fc, \bbF = \{\Fc_t\}_{t\ge 0}, \bbP)$ be a given complete probability space, $W(s)$ be a $d$-dimensional Brownian motion on $(\Omega, \Fc, \bbF, \bbP)$ where $s$ denotes the time. 
For a fixed initial state $\bfx$ at some time $0<t<T<\infty$, we seek to control the randomly perturbed dynamical system  \begin{equation}\label{eq:charOC_stoc}
	\begin{cases}
    		\du \bfz_{t,\bfx}(s) = f(s,\bfz_{t,\bfx}(s), \bfu_{t,\bfx}(s, \bfz_{t,\bfx}(s)))\du s + \sigma(s, \bfz_{t,\bfx}(s)) \du W(s), s \in [t,T], \\ 
    		\bfz_{t,\bfx}(t)= \bfx.
	\end{cases}
\end{equation}
Here, $\bfz_{t,\bfx} : [t,T]\ \to \R^d$ describes the state and $\bfu_{t,\bfx}: [t,T]\times  \R^d \to U$ describes the control of the system, 
the function $\sigma: [t,T]\times \R^d \to \R^{d\times d}$ represents the diffusion coefficient, and $f:[t,T]\times \R^d \times U\to \R^d$ represents the drift of the system. 
We assume that the set of admissible controls $U\subset \R^k$ is closed. 
We seek to minimize the objective functional 
\begin{equation}\label{eq:Joc_stoc}
	J_{t,\bfx}(\bfu_{t,\bfx}) = \bbE \left\{G\big(\bfz_{t,\bfx}(T)\big) + \int_{t}^T L\big(s,\bfz_{t,\bfx}(s), \bfu_{t,\bfx}(s, \bfz_{t,\bfx}(s))\big) \, \du s\right\},
\end{equation}
which is comprised of the running cost $L:[t,T]\times \R^d \times U\to \R$ and the terminal cost $G:\R^d\to \R$.
Here, the expectation is taken with respect to perturbances of the dynamics~\cref{eq:charOC_stoc} that is described by the probability measure $\bbP$. We assume sufficient regularity conditions on $f$, $\sigma$, $G$, and $L$, see \cite[Chapter~2]{yong_zhou-1999} for a list of assumptions.

The value function assigns the optimal cost-to-go to any initial state, that is, 
\begin{equation}\label{eq:OC_stoc}
	\Phi(t,\bfx)=\inf_{\bfu_{t,\bfx}}  J_{t,\bfx}(\bfu_{t,\bfx}),
\end{equation}
and a solution $\bfu_{t,\bfx}^*$ to \cref{eq:OC_stoc} incurring this minimum value is called an optimal control.

The \textit{generalized Hamiltonian}, $H :[t,T] \times \R^d \times \R^d \times \R^{d\times d} \to \R \cup \{\infty\}$, is a key ingredient for the SOC theory in the following sections and the backbone of our numerical scheme. For the problem defined in~\cref{eq:charOC_stoc} and \cref{eq:Joc_stoc} it reads
\begin{equation}\label{eq:H_stoc}
	\begin{split}
		H(s,\bfz,\bfp,\bfM)=\sup_{\bfu\in U} \Hc(s,\bfz,\bfp,\bfM,\bfu),
	\end{split}
\end{equation}
where $\bfp$ and $\bfM$ are called adjoint variables and
\[ \Hc(s,\bfz,\bfp,\bfM,\bfu) =  \frac{1}{2}\text{tr}\left(\sigma(s,\bfz)\bfM\right) + \bfp \cdot f(s,\bfz,\bfu)-L(s,\bfz,\bfu).  \]
We assume that there exists a unique minimizer of the Hamiltonian \cref{eq:H_stoc}.

To make it notationally convenient, in the rest of the paper, we drop the second argument for the controls and denote controls by $\bfu_{t,\bfx}(s)$.

%%%%%%%%%%%%%%%%%%%%%%%%%%%%
\subsection{Stochastic Pontryagin Maximum Principle}
\label{ssec:PMP}
%%%%%%%%%%%%%%%%%%%%%%%%%%%%
The PMP provides first-order necessary conditions for the SOC problem and also states that the optimal control $\bfu_{t,\bfx}^*$ must satisfy an (extended) Hamiltonian system along the optimal state and adjoint trajectory. 
This is made precise by the following result from \cite[Theorem~3.2, Chapter~3]{yong_zhou-1999}. 
\begin{theorem}\cite[Theorem~3.2, Chapter~3]{yong_zhou-1999}\label{thm:PMP}
	Assume that $(\bfz^*_{t,\bfx},\bfu^*_{t,\bfx})$ is an optimal pair that solves~\cref{eq:charOC_stoc} and \cref{eq:Joc_stoc}. Then there exist adjoint states $\bfp_{t,\bfx} \colon [t,T]\to \R^d$ and $\bfM_{t,\bfx} \colon [t,T]\to \R^{d \times d}$ satisfying the adjoint equation
	\begin{equation}\label{eq:adjoint}
		\begin{cases}
			\du \bfp_{t,\bfx}(s)= & \bfM_{t,\bfx}(s) \du W(s) - \nabla_{\bfz} \mathcal{H} \big( s,\bfz^*_{t,\bfx}(s), \bfp_{t,\bfx}(s), \bfM_{t,\bfx}(s),\bfu^*_{t,\bfx}(s)) \du s \\
			\bfp_{t,\bfx}(T)= & - \nabla_{\bfz} G \big(\bfz_{t,\bfx}^*(T) \big),
		\end{cases}
	\end{equation}
	where $s \in [t,T]$ and the optimal control satisfies 
	\begin{equation}\label{eq:PMP}
		\begin{split}
			\bfu^*_{t,\bfx}(s) = \argmax_{\bfu \in U}~ \mathcal{H} \big( s,\bfz_{t,\bfx}^*(s), \bfp_{t,\bfx}(s), \bfM_{t,\bfx}(s),\bfu(s) \big)
		\end{split}
	\end{equation}
	for almost all $s\in [t,T]$, $\bbP$-almost surely.
\end{theorem}
We note that the optimal control defined in \cref{eq:PMP} only depends on the adjoint variable $\bfp_{t,\bfx}$ but not on $\bfM_{t,\bfx}$ since $\sigma(\cdot,\cdot)$ does not depend on the control.

We further assume that there exists a unique continuous closed-form solution to \cref{eq:PMP}.
Although not demonstrated in this work, this assumption can be weakened to include implicitly defined functions as long as they can be obtained efficiently; this allows, for example, modeling more general convex running costs.

We note that when the control satisfies~\cref{eq:PMP}, the dynamics in \cref{eq:charOC_stoc} is equal to
\begin{align}\label{eq:charOC_stoc_red}
	\left\{ \begin{aligned}
		\du \bfz_{t,\bfx}^*(s)& = \nabla_{\bfp} \Hc \big( s,\bfz_{t,\bfx}^*(s), \bfp_{t,\bfx}(s), \bfM_{t,\bfx}(s), \bfu^*_{t,\bfx}(s) \big) \du s	+ \sigma( s,\bfz_{t,\bfx}^*(s)) \du W(s), \\ %\; t \leq  s \leq  T, \\ %
		\bfz_{t,\bfx}^*(t) &= \bfx. 
	\end{aligned} 
	\right.
\end{align}

The system of equations \cref{eq:charOC_stoc_red,eq:adjoint,eq:PMP} is called the stochastic Hamiltonian system, where the maximum condition \cref{eq:PMP} corresponds to the variational inequality for the control.

Finding a tuple
$(\bfz^*_{t,\bfx},\bfu^*_{t,\bfx},\bfp_{t,\bfx},\bfM_{t,\bfx})$ that satisfies the PMP can be extremely difficult. 
However, when the value function $\Phi$ is differentiable, $(\bfp_{t,\bfx},\bfM_{t,\bfx})$ satisfying \cref{eq:adjoint} can be obtained from $\Phi$; this is formalized in the following theorem that, with weaker assumptions, can be found in  \cite[Chapter~5]{yong_zhou-1999}.
\begin{theorem}\cite[Chapter~5]{yong_zhou-1999}\label{thm:p=gradPhi}
	Assume that $\bfu^*_{t,\bfx}$ is an optimal control and $\Phi\in C^{1,3}([t,T]\times\R^d)$. Then
	\begin{equation}\label{eq:gradPhi}
		\bfp_{t,\bfx}(s) = -\nabla_{\bfz} \Phi \big(s,\bfz_{t,\bfx}^*(s)\big) 
 		\; \text{ and } \; 
 		\bfM_{t,\bfx}(s))= - \sigma(s,\bfz_{t,\bfx}^*(s))^\top \nabla^2_{\bfz} \Phi \left(s,\bfz_{t,\bfx}^*(s)\right)
	\end{equation}
	solve \cref{eq:adjoint}.
\end{theorem}
\Cref{thm:p=gradPhi} along with \cref{eq:PMP} collectively serve to express the optimal control as
\begin{equation}\label{eq:opt_control_nec}
	\begin{split}
		\bfu^*_{t,\bfx}(s)=\bfu_{t,\bfx}^*\left(s,\bfz^*_{t,\bfx}(s),-\nabla_{\bfz} \Phi \big(s,\bfz_{t,\bfx}^*(s)\big)\right).
	\end{split}
\end{equation}
This relation along with \cref{eq:charOC_stoc_red} is one of the key ingredients of our numerical solution approach.
\Cref{eq:opt_control_nec} characterizes optimal control in a feedback or closed-loop form, which is of utmost importance in many real-life applications. This is because optimal control can be quickly computed at any given point in time and space  when the value function is known and its gradient is readily available. Thereby, avoiding re-computation of optimal controls for cases when multiple evaluations are needed for different times or states.

%%%%%%%%%%%%%%%%%%%%%%%%%%%%
\subsection{Hamilton-Jacobi-Bellman Equation}
\label{ssec:HJB}
%%%%%%%%%%%%%%%%%%%%%%%%%%%%
To help approximate  the value function $\Phi$, we also use the fact that  $\Phi$ satisfies the Hamilton-Jacobi-Bellman (HJB) equation, which is a result of the Dynamic Programming (DP) method or Bellman's principle. 
We state the following result taken from \cite{yong_zhou-1999} under suitable assumptions, see also \cite[Remark~3.4.4, Theorem~3.5.2]{pham-2009}.
\begin{theorem}\cite[Propositon~3.5, Chapter~4]{yong_zhou-1999}
	Assume that the value function $\Phi \in C^{1,2}([t,T] \times \R^d)$. Then $\Phi$ satisfies the HJB equation
	\begin{alignat}{2}\label{eq:HJB_stoc}
		\begin{split}
			-&\partial_s \Phi(s,\bfz)  +  H\big(s,\bfx,-\nabla_{\bfz} \Phi(s,\bfz),  -\sigma(s,\bfz)^\top\nabla_{\bfz}^2 \Phi(s,\bfz)\big)=0,\; \forall (s,\bfz)\in [t,T)\times \R^d,\\
			&\Phi(T,\bfz) = G(\bfz).
		\end{split} 
	\end{alignat}
\end{theorem}
The smoothness of $\Phi$ can be relaxed to continuity in the weaker sense of viscosity solutions~\cite[Section~5, Chapter~4]{yong_zhou-1999}.

Using the definition of the Hamiltonian in \cref{eq:HJB_stoc} we get the HJB equation as the following second-order parabolic PDE:
\begin{alignat}{2}\label{eq:HJB_pde}
	\begin{cases}
		-\partial_s \Phi(s,\bfz) -\dfrac{1}{2}\text{tr}(\sigma(s,\bfz)\sigma(s,\bfz)^\top\nabla^2_{\bfz} \Phi(s,\bfz)) - \nabla_{\bfz} \Phi(s,\bfz)\cdot f(s,\bfz,\bfu^*)\\
		\phantom{-\partial_s \Phi(t,\bfx) -\dfrac{\sigma^2}{2}\Delta \Phi(t,\bfx) + \nabla \Phi(t,\bfx)-}-L(s,\bfz,\bfu^*)=0,\; \forall \,(s,\bfz)\in [t,T)\times \R^d,\\
		\Phi(T,\bfz) = G(\bfz(T)).
	\end{cases}
\end{alignat}
In addition, by the envelope theorem, it follows that $\nabla_{\bfp} \Hc = \nabla_{\bfp} H$ and $\nabla_{\bfM} \Hc = \nabla_{\bfM} H$. This simplifies the computation of optimal trajectories, which can now be expressed via the value function as:
\begin{align}\label{eq:HJB_FSDE}
	\left\{ \begin{aligned}
		\du \bfz_{t,\bfx}^*(s)& = \nabla_{\bfp} H \big( s,\bfz_{t,\bfx}^*(s), -\nabla_{\bfz} \Phi \big(s,\bfz_{t,\bfx}^*(s)\big), -\sigma(s,\bfz_{t,\bfx}^*(s))^\top \nabla^2_{\bfz} \Phi \big(s,\bfz_{t,\bfx}^*(s)\big) \big) \du s \\  & \phantom{\nabla_{\bfp} H \big( s,\bfz_{t,\bfx}^*(s), -\nabla_{\bfz} \Phi \big(s,\bfz_{t,\bfx}^*(s)\big)\big)\sigma(s,\bfz_{t,\bfx}^*(s))} + \sigma(s,\bfz_{t,\bfx}^*(s))  \du W(s), \\ %
		\bfz_{t,\bfx}^*(t) &= \bfx.
	\end{aligned} 
	\right.
\end{align}
These modified dynamics do not explicitly involve the control, which reduces the problem solely to the state variables. 
Note that in \cref{eq:HJB_FSDE} the dependency of $H$ on the Hessian of the value function could be omitted because the volatility term does not depend on the control.
This idea was proposed for deterministic optimal control problems in \cite{onken2020neural}.

%%%%%%%%%%%%%%%%%%%%%%%%%%%%
\subsection{FBSDE Formulation}
%%%%%%%%%%%%%%%%%%%%%%%%%%%%
One way to avoid the need for a spatial discretization of the HJB equation~\cref{eq:HJB_pde} is to use a non-linear version of Feynman-Kac formula and obtain an equivalent system of stochastic differential equations. 
This idea has been applied to a variety of nonlinear parabolic/elliptic PDEs; see,  for example, \cite{Antonelli_1993,Ma_Protter_Yong_1994,Cheridito_Touzi_2007,Pardoux_Tang_1999,pham-2009}.

Our FBSDE system uses ~\cref{eq:HJB_FSDE} as the forward system to sample trajectories. 
Along those trajectories, we note that the solution to the HJB equation \cref{eq:HJB_pde} must satisfy the backward SDE following the Feynman-Kac formulae, see \cite[Chapter~7]{yong_zhou-1999}
\begin{alignat}{3} \label{eq:HJB_BSDE}
	\begin{split} 
		\begin{cases}
			\du \Phi(s,\bfz(s))& = \nabla_{\bfz} \Phi(s,\bfz(s))^\top \sigma(s,\bfz(s)) \,\du W(s)
			- L(s,\bfz(s), \bfu^*(s))\du s, \quad \\
			\Phi(T,\bfz(T))& = G(\bfz(T)).
		\end{cases}	
	\end{split} 
\end{alignat}

It is important to stress that our choice of the forward system \cref{eq:HJB_FSDE} is the key difference from other existing works.
For example, ~\cite{Han_BSDE} and \cite{raissi2018forward} use the standard Brownian motion. 
While both of these choices lead to a valid FBSDE system for \cref{eq:HJB_pde}, we advocate for including the control in the dynamics as motivated by stochastic PMP \cref{eq:charOC_stoc_red}. As our numerical experiments demonstrate, focusing the sampling along optimal trajectories can lead to more accurate and efficient value function approximations.

%%%%%%%%%%%%%%%%%%%%%%%%%%%%%%%%%%%%%%%
\section{Neural Network Approach}
\label{s:NN-approach}
%%%%%%%%%%%%%%%%%%%%%%%%%%%%%%%%%%%%%%%
In this section, we present a neural network framework for approximating the value function of the stochastic optimal control problem defined by the objective functional \cref{eq:Joc_stoc} and dynamics \cref{eq:charOC_stoc}.
The theoretical foundation of our framework is given by the PMP, FBSDE, and Dynamic Programming as presented in the previous section. 
The key idea is to approximate the value function $\Phi$ in \cref{eq:OC_stoc} by a neural network and compute the control using the feedback form \cref{eq:opt_control_nec}. 
What distinguishes our framework from similar approaches such as~\cite{Han_BSDE,raissi2018forward} is the use of the feedback form to guide the sampling during training.
Thereby we seek to learn to explore the relevant part of the state space.
We also derive and experiment with various loss functions that are based on the control objective \cref{eq:Joc_stoc}, the BSDE \cref{eq:HJB_BSDE}, and the  HJB~\cref{eq:HJB_pde}.

\subsection{Neural Network Approximation}
The first building block of our framework is to parameterize the value function using a neural network.
Since finding an effective network architecture for any learning task is both crucial and an open research topic, we treat this as a modular component.
Our framework can be used with any  scalar-valued neural network that takes inputs in $\R^{d+1}$ as long as it is twice differentiable with respect to its last $d$ inputs; this is to allow computations of $\nabla \Phi$.

Among the networks we use in our numerical experiments is the multi-layer perceptron (MLP) model used in~\cite{raissi2018forward}.
As an alternative, which also satisfies the regularity needed, we propose the residual network also used for deterministic control in~\cite{onken2020neural}. The network is given by
\begin{equation}
\label{eq:NNArchitecture}
	% \begin{split}
	\Phi(\bfy ; \bfth) = \bfw^\top \mathcal{NN}(\bfy;\bfth_\mathcal{NN}) 
	+ \frac{1}{2} \bfy^\top (\bfA^\top 
	\bfA)\bfy + \bfb^\top \bfy + c, 
	% \end{split}
\end{equation}
with trainable weights $\bfth = (\bfw, \bfth_\mathcal{NN}, \bfA, \bfb, c)$.
Here the inputs $\bfy=(s,\bfz)\in \R^{d+1}$ correspond to time-space,  $\mathcal{NN}(\bfy;\bfth_\mathcal{NN}) \colon \R^{d+1} \to \R^m$ is a neural network,
	and $\bfth$ contains the trainable weights: $\bfw\,\,{\in}\,\,\R^m$, $\bfth_\mathcal{NN}\,\,{\in}\,\,\R^p$
	, $\bfA\,\,{\in}\,\,\R^{\gamma \times (d+1)}$, $\bfb\,\,{\in}\,\,\R^{d+1}$, $c\,{\in}\,\R$, where rank $\gamma{=}\min (10,d+1)$ limits the number of parameters in $\bfA^\top \bfA$. 
	Here, $\bfA$, $\bfb$, and $c$ model quadratic potentials, that is, linear dynamics; $\mathcal{NN}$ models nonlinear dynamics. For certain experiments, we may choose to omit the quadratic potential terms $\bfA$, $\bfb$ and $c$ for comparison or simplicity reasons.
 
	In our experiments, for $\mathcal{NN}$, we either use a MLP \cite{Goodfellow-et-al-2016}
	\begin{equation} \label{eq:MLP}
	    \begin{split}
	    \bfa_0 & = \rm act(\bfK_0 \bfy + \bfb_0) \\ 
	    \bfa_{i+1} & = \rm act(\bfK_{i+1} \bfa_i + \bfb_{i+1}),\quad 0\le i \le M-2 \\
	    \mathcal{NN}(\bfy;\bfth_\mathcal{NN}) & =  \rm act(\bfK_M \bfa_{M-1} + \bfb_M),
	    \end{split}
	\end{equation}
	or a residual neural network (ResNet)~\cite{He_2016_CVPR}	
	\begin{equation} \label{eq:ResNet}
	    \begin{split}
	    \bfa_0 & = \rm act(\bfK_0 \bfy + \bfb_0) \\ 
	    \bfa_{i+1} & = \bfa_i + \rm act(\bfK_{i+1} \bfa_i + \bfb_{i+1}),\quad 0\le i \le M-2 \\
	    \mathcal{NN}(\bfy;\bfth_\mathcal{NN}) & = \bfa_{M-1} + \rm act(\bfK_M \bfa_{M-1} + \bfb_M), 
	    \end{split}
	\end{equation}	
	with neural network weights $\bfth_\mathcal{NN}{=}(\bfK_0,\hdots,\bfK_M, \bfb_0,\hdots,\bfb_M)$ where $\bfb_i\in \R^{m}\; \forall i$, $\bfK_0 \in \R^{m \times (d+1)}$, and $\{ \bfK_1,\hdots,\bfK_M \}\in \R^{m \times m}$ with $M$ being the depth of the network.
	The choice of the element-wise nonlinearity $\rm act(\cdot)$ is discussed in the respective  experiments. 

\subsection{Training Problem}
Ideally, we would choose $\bftheta$ such that $\Phi(s,\bfz;\bftheta)$ is equal to the value function of the control problem  globally, that is, for all $(s,\bfz) \in [t,T]\times \R^d$.
Since this is known to be cursed by the dimensionality for reasonable problem sizes, we resort to a semi-global approach, which enforces this property at randomly sampled points in the space-time domain.

To generate samples, we first obtain initial states $\bfx\sim \rho$ from some (possibly Dirac)  distribution $\rho$ and then use an Euler Maruyama scheme with  $N+1$ equidistant time points $s_0, \ldots, s_N$ and step size $\du s = (T-t)/N$.
This yields a state trajectory starting at $\bfz_0=\bfx$ via
\begin{equation}\label{eq:sampling}
% \begin{align}
    \bfz_{i+1} = \bfz_i + f(s_i,\bfz_i,\bfu_i) \du s + \sigma(s_i,\bfz_i) \du \bfW_i, \quad i=0,\ldots,N-1
% \end{align}    
\end{equation}
where $\du\bfW_i \sim \mathcal{N}({\bf 0},\du s \cdot \bfI_d)$,  and $\bfu_i = \bfu_{t,x}^*(s_i,\bfz_i)$ is the optimal control obtained from the feedback, that is, form~\cref{eq:PMP}
\begin{equation*}
    \bfu_i^* \in \argmax_{\bfu \in U} \Hc\left(s_i,\bfz_i,-\nabla \Phi(s_i,\bfz_i;\bftheta),-\sigma(s_i,\bfz_i)^\top\nabla^2 \Phi(s_i,\bfz_i;\bftheta), \bfu\right).
\end{equation*}
A few comments are in place. First, it is important to note that due to the feedback form, the sampled trajectories depend on the parameters of the value function. 
Second, the addition of this drift term, motivated by control theory, is the key difference to neural network solvers for the more general class of semi-linear elliptic PDEs~\cite{Han_BSDE,raissi2018forward}.
Third, the drift term can also be motivated by the fact that for $\sigma \to 0$, the trajectories above approximate the characteristic curves of the non-viscous HJB equation; thereby our SOC approach coincides with that for deterministic OC in~\cite{onken2020neural}.

To further simplify the notations, we omit the subscript $\bfz$ in $\nabla_{\bfz} \Phi$ and $\nabla^2_{\bfz} \Phi$ for the rest of the paper. Furthermore, we collect the states, control, and noise along the discrete trajectories in~\cref{eq:sampling} column-wise in the matrices
\begin{equation*}
    \bfZ \in \R^{d\times N}, \quad \bfU \in \R^{k\times N}, \quad \mathbf{dW} \in \R^{d\times N}.
\end{equation*}

To learn the parameters of the neural networks in an unsupervised way (that is,  assuming neither analytic values of $\Phi$ nor optimal control trajectories), we approximately solve the minimization problem
\begin{equation}\label{eq:full_opt}
\begin{split}
    	\min_{\bfth} \E_{\bfx \sim \rho}  & \left\{\E_{\bfZ,\bfU,\mathbf{dW} | \bfx} \left\{
      \beta_1 P_{\rm BSDE}^{p}(\bfZ,\bfU,\mathbf{dW}) + \beta_2 P_{\rm HJB}^{p}(\bfZ) + \beta_3 J(\bfZ,\bfU) \right.\right.\\
      &\left.\left.
       +\beta_4 |G(\bfz_N) - \Phi(s_N, \bfz_N; \bftheta)|^p +  \beta_5 |\nabla G(\bfz_N) - \nabla\Phi(s_N,\bfz_N; \bftheta)|^p
      \right\} \right\},
      \end{split}
\end{equation}
where the  terms in the objective function consist of penalty functions for violations of the BSDE system and the HJB PDE, the control objective, and penalty terms for the terminal condition, respectively, and are defined below. The exponent $p\in\{1,2\}$ allows one to choose between different norms for the loss function. In our numerical examples, we use $p = 1$ as it gives much faster convergence. The relative influence of each term is controlled by the components of $\beta \in \R^5_{+}$.
Different choices of $\beta$ allow us to experiment with different learning approaches; for example, setting $\beta_1 = \beta_4 = \beta_5=1$ and $\beta_2=\beta_3=0$ provides the same loss function as in~\cite{raissi2018forward} while $\beta_1=0$ and $\beta_i>0, \; i\in \{2,3,4,5\}$ gives the loss function used for deterministic OC problems in~\cite{onken2020neural}.

We penalize the violation of the BSDE~\cref{eq:HJB_BSDE} via
\begin{equation}\label{eq:p_bsde}
    P_{\rm BSDE}(\bfZ,\bfU,\mathbf{dW}) = \sum_{i=0}^{N-1} |\Phi_{i+1}(\bftheta)-\Phi_{i}(\bftheta) + L(s_i,\bfz_i,\bfu_i)\du s -  \nabla\Phi_i(\bftheta)^\top \sigma(s_i,\bfz_i)  \du\bfW_i |
\end{equation}
where we use the abbreviations $\Phi_i(\bftheta) := \Phi(s_i,\bfz_i;\bftheta)$ and $\nabla\Phi_i(\bftheta) := \nabla\Phi(s_i,\bfz_i;\bftheta)$.
Similarly, the HJB penalty term reads
\begin{equation}\label{eq:p_hjb}
% \begin{split}
     P_{\rm HJB}(\bfZ) = \du s \sum_{i=1}^{N} |H(s_i,\bfz_i,-\nabla \Phi_i(\theta), -\sigma(s_i,\bfz_i)^\top \nabla^2 \Phi_i(\theta)) - \partial_s \Phi_i(\bftheta)|,
% \end{split}
\end{equation}
where $\nabla^2 \Phi_i(\theta) := \nabla^2\Phi(s_i,\bfz_i;\theta)$,  $\partial_s \Phi_i(\theta) := \partial_s \Phi(s_i,\bfz_i;\theta)$. 
Finally, we approximate the objective functional via
\begin{equation*}
    J(\bfZ,\bfU) = G(\bfz_N) + \du s \sum_{i=1}^N L(s_i,\bfz_i,\bfu_i).
\end{equation*}

In principle, any stochastic approximation approach can be used to approximately solve the above optimization problem.
Here, we use Adam~\cite{kingma2014adam} and sample a minibatch of trajectories originating in i.i.d. samples from $\rho$.

\subsection{Implementation}
We implement and test our proposed approach in two software environments.

To obtain a direct comparison with~\cite{raissi2018forward} we modify the FBSNN code accompanying the paper.
To this end, we created a publicly available fork at \url{https://github.com/EmoryMLIP/FBSNNs}. 
Our two main modifications are adding the proposed drift to the forward dynamics and adding the control objective in the training loss. 
Other parameters, including the choice of neural network model, are kept unchanged.

In order to further simplify the experimentation, we also implement our own PyTorch code available at \url{https://github.com/EmoryMLIP/NeuralSOC}. Our implementation contains all loss terms in \cref{eq:full_opt}. 
We implement both sampling techniques: pure random walk and the proposed one informed by PMP.
This facilitates comparisons of  our approach with other available methods and simplifies developing new examples. 

We tested most of our examples using either Intel Xeon E5-4627 CPU or Nvidia P100 GPU.

%%%%%%%%%%%%%%%%%%%%%%%%%%%%%%%%%
\section{Numerical Experiments}
\label{sec:num_exp}
%%%%%%%%%%%%%%%%%%%%%%%%%%%%%%%%

We test the efficacy of our proposed algorithm on several different SOC problems. 
In \cref{2D problem}, we  introduce a two-dimensional trajectory planning problem to visualize the difference between purely random exploration and our proposed sampling scheme.  To illustrate the accuracy of the learned value function, we compare it with the value function obtained by solving the corresponding HJB PDE using  a finite element method (FEM). The goal of this experiment is to compare the accuracy of the neural network and FEM approximation. 
In \cref{s:100d-exp}, we compare our approach to those in \cite{Han-2017,Han_BSDE} using a  100-dimensional benchmark problem. For the original version of this problem, our method shows faster initial convergence and time-to-solution with comparable accuracy.
We modify the terminal cost of this problem to further highlight the importance of the feedback form in the sampling  (see \cref{s:100d-exp-mod}). Lastly in \cref{sec:quadcop} we also test our method on a 12-dimensional problem with nonlinear dynamics, showing that our method generates relatively accurate solutions under complex dynamics.

\subsection{2D Trajectory Planning Problem}
\label{2D problem}
To visualize the behavior of our PMP-based sampling approach, we consider a two-dimensional test problem. 

The problem consists of planning an optimal trajectory from the initial state $\bfx \sim \rho = \mathcal{N}((-1.5,-1.5)^\top, 0.4 \cdot \bfI_2)$ to the target $\bfx_{\rm target} = (1.5,1.5)^\top$. 
To make the problem interesting, a hill is placed at the origin, denoted by $Q(\bfz)$, which adds height-dependent cost for traveling around that region. In our experiments, $Q(\bfz)$ is defined by a two-dimensional Gaussian density with mean zero and covariance of $0.4 \cdot \bfI_2$ scaled by a factor of 50.

\begin{figure}[t]
\centering
\includegraphics[width = 0.35\textwidth]{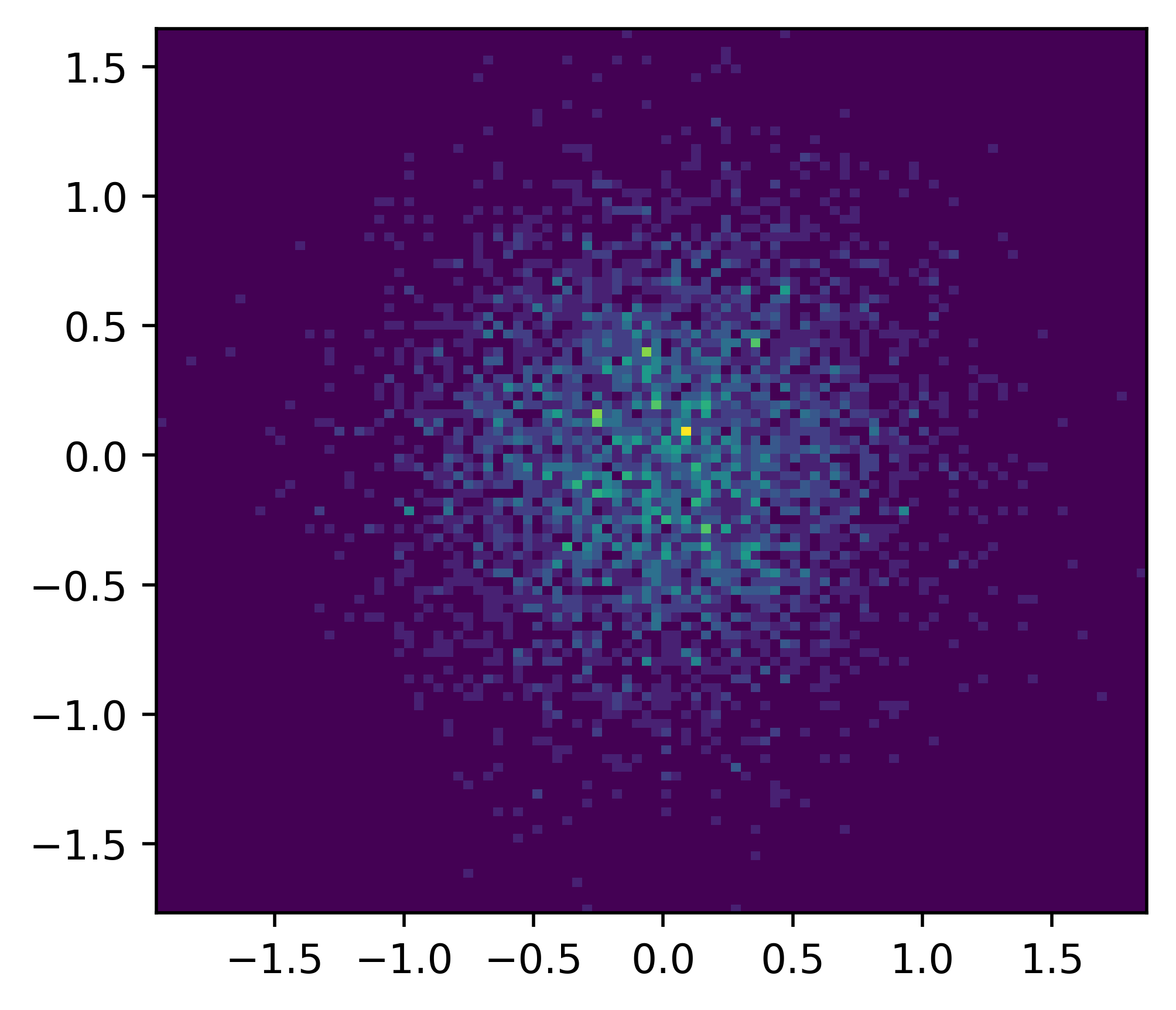}%
\includegraphics[width = 0.35\linewidth]{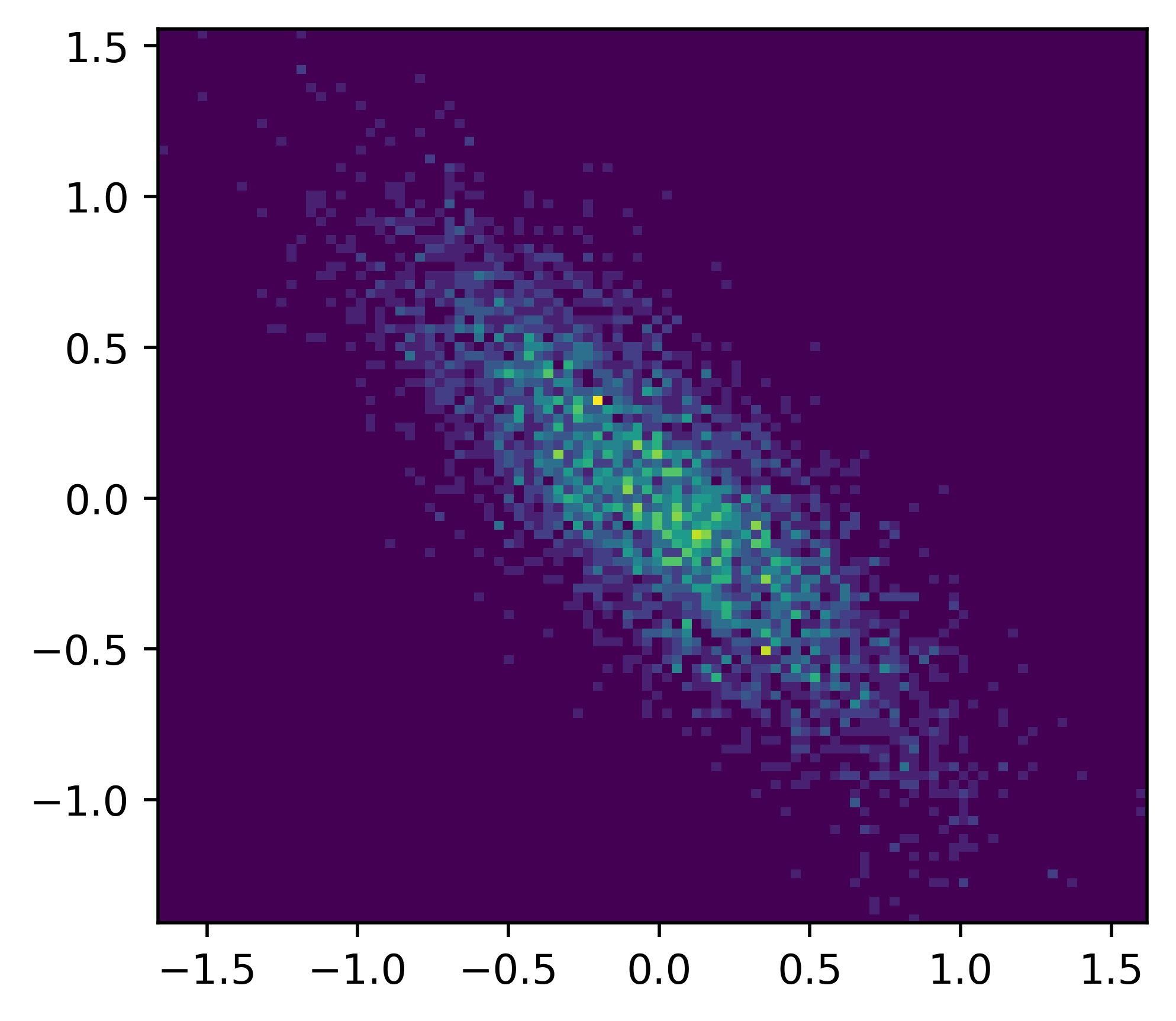}%
\caption{Action of $\sigma$ in \eqref{eq:2d_dyn} on standard Gaussian distribution (Left) warps it diagonally (Right). This would affect the solution of the problem.}
\label{fig:sigma matrix}
\end{figure}

The dynamics for the problem read 
\begin{equation}\label{eq:2d_dyn}
	f(s,\bfz,\bfu) = \bfu \quad \text{and} \quad \sigma = \begin{bmatrix} 0.2 & -0.4 \\ -0.4 & 0.2 \end{bmatrix}.
\end{equation}
The choice of non-scalar $\sigma$ adds to the complexity of the problem by changing the behavior of the standard Brownian motion, see \cref{fig:sigma matrix}.

The running cost and terminal cost of the problem are given, respectively, by
\begin{equation}\label{eq:2d_LG}
	L(s,\bfz,\bfu) = \frac{1}{2} \| \bfu \|^2 + Q(\bfz) \quad \text{and} \quad  G(\bfz) = 50 \cdot \| \bfz - \bfx_{\rm target} \|^2.
\end{equation}

The corresponding HJB equation is
\begin{subequations}\label{eq: 2d prob}
\begin{equation}
	\partial_s\Phi (s, \bfz) + \frac{1}{2} \text{tr} (\sigma \sigma^\top \nabla^2 \Phi(s, \bfz)) - \frac{1}{2}\| \nabla \Phi(s, \bfz) \|^2 + Q(\bfz) = 0.
\end{equation}
with terminal condition 
\begin{equation}
    \Phi(T, \bfz) = G(\bfz).
\end{equation}
\end{subequations}

\subsubsection{Finite Element Method}
Since it is not obvious how to solve the HJB equation \cref{eq: 2d prob} analytically, we approximately solve it using a finite element method (FEM) to obtain a baseline for this problem.

We approximate the value function by solving the HJB PDE \cref{eq: 2d prob} on the domain $\Omega = [-3,3]\times [-3,3]$  with homogeneous Neumann boundary conditions,
\begin{alignat*}{3}
    	&\frac{\partial \Phi}{\partial \hat{\bfn}}(s,\bfz)  = 0, \; & \text{on } \partial \Omega,\; \forall s<T, 
\end{alignat*}
where $\hat{\bfn}$ denotes the unit normal vector.
Since the diffusion coefficient $\sigma$ is independent of time and space, $\text{tr} (\sigma \sigma^\top \nabla^2 \Phi(s, \bfz)) = \text{div}(\sigma \sigma^\top \nabla \Phi(s, \bfz))$, which we use to derive a weak form of the PDE. 
Using the implicit Euler discretization in time on a partition of $[0,T]$ into $N$ sub-intervals with uniform step size, $\du s$, yields
\[ \frac{\Phi_{n+1} - \Phi_n}{ \du s } + \frac{1}{2} \text{div} (\sigma \sigma^\top \nabla \Phi_n) - \frac{1}{2}\| \nabla \Phi_n \|^2 + Q = 0, \quad n=N, N-1, \ldots, 0, \]
where $\Phi_n$ denotes the approximated solution $\Phi(t_n,\cdot)$, at $t_n=n\cdot \du s$ and $\Phi_{N+1} = G(\cdot)$.
Then, using  Green's formula, the weak problem at the $n$-th time step consists of finding $\Phi_{n} \in H^1(\Omega)$  such that 
\begin{align*}
   &\int_{\Omega} (\Phi_{n+1} - \Phi_n) v \du \bfz  - \du s\frac{1}{2} \int_{\Omega} \sigma \sigma^\top \nabla\Phi_n \cdot \nabla v \du \bfz + \du s \int_{\Omega} \left(Q- \frac{1}{2} \| \nabla \Phi_n) \|^2\right) v \du \bfz = 0,
\end{align*}
for all test functions $ v\in H^1(\Omega)$.
Here, $H^1(\Omega)$ denotes the Hilbert Sobolev space defined by $H^1(\Omega) = \{v\in L^2(\Omega)| \nabla v \in L^2(\Omega)\}$.

To solve the weak problem we use FEniCS \cite{LangtangenLogg2017}, we create a triangular mesh for $\Omega$ and use $\mathcal{P}_1$ Lagrange finite elements to discretize $\Phi$ in space. We discretize $\Omega$ using 150 mesh points in each dimension, summing up to a total of 22,500 degrees of freedom, and use the step-size of $\du s = 0.001$ in time. At each time step, we use Newton's method to solve for $\Phi_n$, with relative error and absolute error tolerance for the solver set to $10^{-6}$ and $10^{-10}$, respectively. We denote the FEM solution by $\Phi_{\rm FEM}$.

\begin{figure}[t]
\centering
\includegraphics[width = 0.95\textwidth]{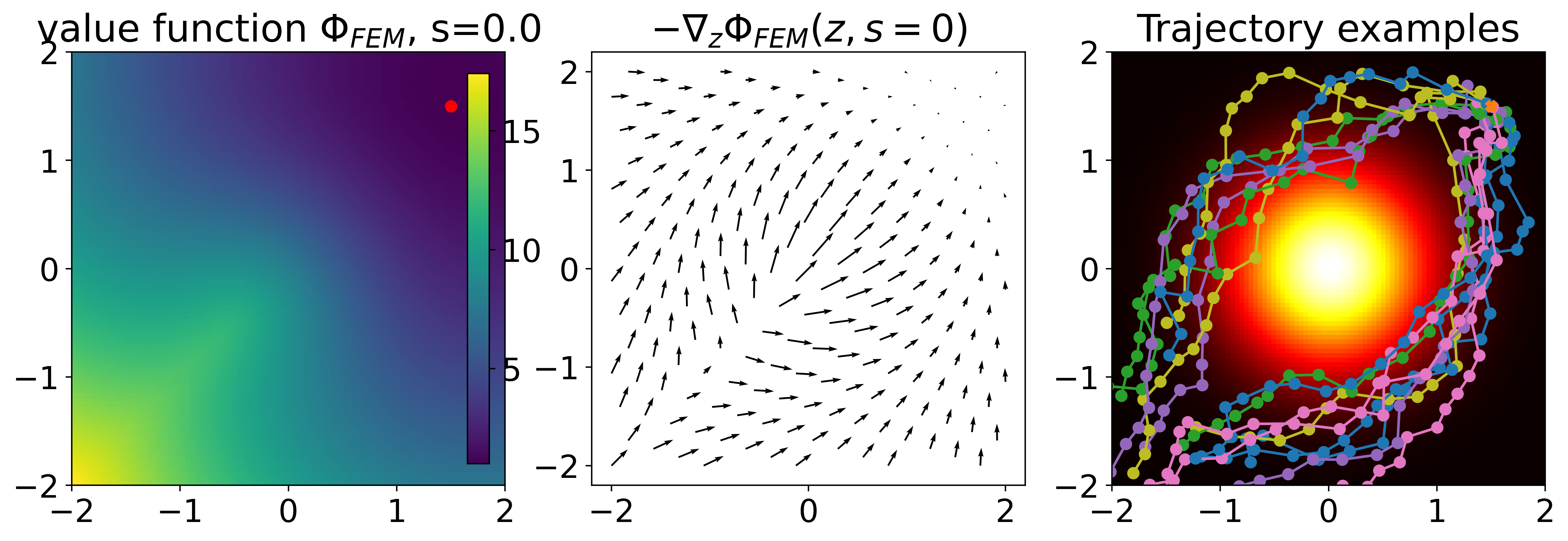}%
\caption{Results of the two-dimensional test problem. Left: Value function approximation $\Phi_{\rm FEM}(0,\cdot)$. Middle: Quiver plot of optimal controls  at $s=0$. Right: Trajectories generated from randomly chosen initial states.}
\label{fig:2d_compile}
\end{figure}

In \cref{fig:2d_compile} we plot the solution $\Phi_{\rm FEM}$ as well as the optimal control policy at initial time $s = 0$  obtained via the feedback form. We also present trajectories originating from some randomly chosen initial states following the optimal policy. As expected, the trajectories travel from the initial points to the target while avoiding the obstacle in the center of the domain.

In \cref{tab:Phi_J}, we evaluate the control objective, $J$, for some fixed initial state. We notice that the estimated value matches with $\Phi_{\rm FEM}(0)$, suggesting that the FEM solution is an accurate approximation of the true $\Phi$.
It is also worth pointing out that FEM is sufficient for the 2D parabolic equation we have here since a variational form is explicitly available. However for problems without a easily accessible variational form, one may want to resort to methods in \cite{M2AN_1995__29_1_97_0,Kushner-1990} for true solutions.

\subsubsection{Neural Network Approach}
For the problem defined in \cref{eq:2d_dyn} and \cref{eq:2d_LG}, the forward SDE~\cref{eq:sampling} simplifies to
\begin{equation}
    \label{eq: 2D_FSDE}
    \bfz_{i+1} = \bfz_{i}- \nabla \Phi(s_i,\bfz_{i})\du s + \sigma \du \bfW_i.
\end{equation}

Following our proposed method in ~\cref{eq:NNArchitecture}, we approximate the value function using a three-layer residual neural network with 32  neurons per layer. We do not include the quadratic terms in the network for this experiment since the simpler structure was already sufficient for solving this problem. Overall, the model consists of 1,217 trainable parameters. We choose $\tanh$ as the activation function for all but the final layer of the network, the final layer does not have an activation function. We use the penalty parameters $\beta=(1.0.1.0,1.0,1.0,0.0)$, that is, we enable both the penalty terms, $P_{\rm BSDE}$ and $P_{\rm HJB}$ in \cref{eq:full_opt}. To approximately solve~\cref{eq:full_opt} we use a total of 6,000 steps of the Adam optimizer with a batch size of 64. We start with a learning rate of $0.01$ and divide it by 10 every 1800 iterations. The average cost per iteration is about  $ 0.22$s using an NVIDIA P100 GPU. Note that for the chosen $\sigma$ in \cref{eq:2d_dyn} full Hessian information of the value function, $\nabla^2 \Phi$, is required to calculate $P_{\rm HJB}$. To this end, we use the efficient implementation in the  package hessQuik \cite{Newman2022}. We refer to the neural network approximated solution as $\Phi_{\rm NN}$ in the following sections. We also notice that  sampling the initial states from a slightly larger area during training often helps the robustness of the learned model. 
Given the stochastic nature of the problem and the random initialization of neural network weights, each training sequence can produce a slightly different model. To account for this,  we repeat the training ten times and obtain neural network approximations of the value functions  $\Phi_{\rm NN}^{(j)}$, where $1 \le j \le 10$. 
We compare the resulting models to the FEM solution in the next section.

\begin{figure}[t]
	\begin{subfigure}[b]{1\linewidth}
		\centering
		\includegraphics[width = 1\linewidth]{2D_new_sigma/Sampling_BFSNN.png}
		\caption{Training samples with pure random walk as FSDE  as also used in \cite{Han_BSDE} and \cite{raissi2018forward}. 
		}
		\label{fig:sampling-deepbsde}
	\end{subfigure}
	\newline
	\begin{subfigure}[b]{1\textwidth}
		\centering
		\includegraphics[width = 1\linewidth]{2D_new_sigma/Sampling_ours.png}
		\caption{Training samples with PMP-based drift term.
		}
		\label{fig:sampling-ours}
	\end{subfigure}
	\caption{We visualize training samples of a pure random walk sampler (top row) and our proposed PMP-based sampler (bottom row) for the two-dimensional test problem.
    At six time points (left to right), we visualize the sampled states as two-dimensional histograms. 
    As expected, the pure random walk explores the area around the initial state in all (even suboptimal) directions, while the proposed approach   learns to sample around approximately optimal trajectories.
 }
	\label{fig:sampling}
\end{figure}
To gain more insight into the sampling, we store all states visited during training and plot them as two-dimensional histograms for different time points (left to right) in  \cref{fig:sampling}.
We compare the proposed PMP-based sampling (\cref{fig:sampling-ours}) to the purely noisy dynamics (\cref{fig:sampling-deepbsde}), that is, without drift, as used in~\cite{Han_BSDE, raissi2018forward}.
As expected, the use of purely noisy dynamics leads to the sampling of points only around the initial states in all (even sub-optimal) directions with almost no samples close to the target.
On the other hand, with the use of drift term, the sampled states visit the paths between the initial and target states.

Another way to interpret the histogram plots in \cref{fig:sampling} is by observing the semi-global nature of our neural network approach for SOC problems.
Since the loss function in \cref{eq:full_opt} penalizes the HJB and BSDE losses in a neighborhood of points sampled using the forward SDE, one would expect the trained model to be more reliable in regions that are frequently visited.

\subsubsection{Comparison}
In this subsection, we compare the neural network models $\Phi^{(1)}_{\rm NN}, \ldots, \Phi_{\rm NN}^{(10)}$ and $\Phi_{\rm FEM}$ along the approximately optimal trajectories. Specifically, we randomly sample initial states from $\rho$ and simulate the trajectories using the trained models. We believe that this approach enables a meaningful comparison since the training procedure focuses on those parts of the state space visited by the trajectories and hence the neural networks approximate the value function semi-globally.

\begin{figure}[t]
    \centering
    \includegraphics[scale=0.37]{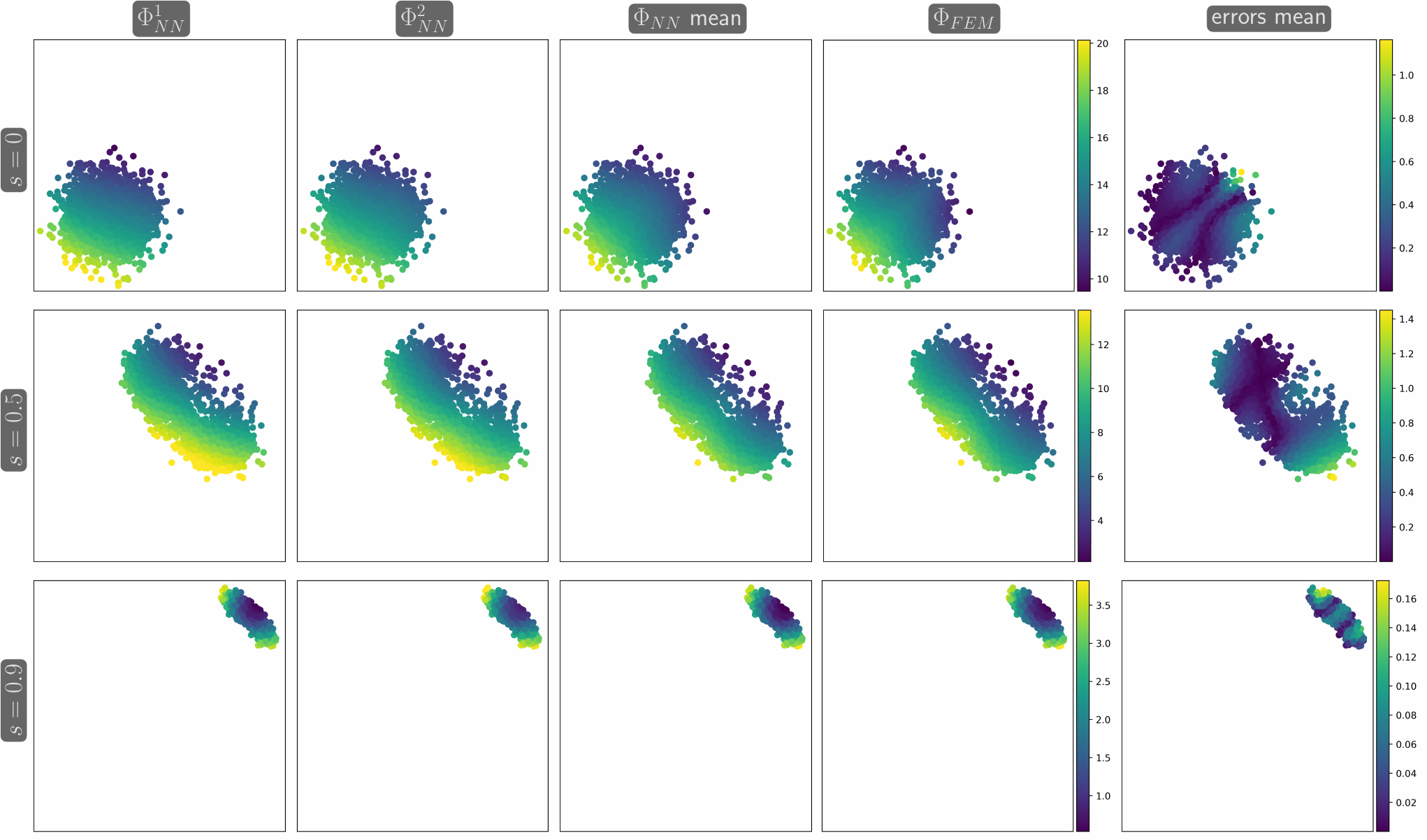}
    \caption{Comparison between learned value function (First 3 columns, including both individual model and model average) and the FEM solution (fourth column) at different time shots $s = 0, 0.5$ and $0.9$. From the errors (fifth column), the neural network solution matches the FEM solution closely over the sampled region.}
    \label{fig:2d_compare}
\end{figure}

For each trained model, we record all the sampled states visited at times $s \in \{0, 0.5, 0.9\}$ while following the corresponding learned policy. We then compare the learned value functions $\Phi_{\rm NN}^j$ with the reference solution $\Phi_{\rm FEM}$ at all these points.
In \cref{fig:2d_compare}, we plot the comparison for times $s \in \{0, 0.5, 0.9\}$ along the rows. The first, second, and fourth columns represent $\Phi_{\rm NN}^{(1)}$, $\Phi_{\rm NN}^{(2})$, and $\Phi_{\rm FEM}$ at the sampled points, respectively. 
We observe that value function estimates look similar. The third column shows the average of the ten learned value functions obtained from the ten training sequences. Lastly, the last column displays the average absolute mean errors between the learned value functions and $\Phi_{\rm FEM}$.

In \cref{tab:error_tab }, we compare mean of the absolute and relative errors between $\Phi_{\rm FEM}$ and $\Phi_{\rm NN}^{(j)}$, $1\le j \le 10$, across the sampled points shown in \cref{fig:2d_compare} for all ten trained models, computed via
\begin{equation}\label{eq:AE}
    {\rm AE}(s) = \frac{1}{n_\text{samples} \times n_\text{models}} \sum_{i = 1}^{n_\text{samples}}\sum_{j = 1}^{n_\text{models}} \left| \Phi_{\rm NN}^{(j)}(s, \bfz_i) - \Phi_{\rm FEM}(s, \bfz_i) \right|
\end{equation}
and
\begin{equation}\label{eq:RE}
   {\rm RE}(s) = \frac{1}{n_\text{samples}\times n_\text{models}} \sum_{i = 1}^{n_\text{samples}}\sum_{j = 1}^{n_\text{models}} \frac{ \left| \Phi_{\rm NN}^{(j)}(s, \bfz_i) - \Phi_{\rm FEM}(s, \bfz_i) \right|}{\left| \Phi_{\rm FEM}(s, \bfz_i) \right|}.
\end{equation}
Our observations indicate that the relative error is smallest at the initial time and increases over time, while the average absolute error remains fairly constant across all states and time intervals. One possible explanation is that we include the control objective, $J$, in the training loss function. Furthermore, the errors across all the trained models display a relatively low standard deviation, indicating that our proposed training scheme is robust to random initialization.

\begin{table}[t]
    \centering
    \begin{tabular}{cccc}
    \hline  \noalign{\smallskip}
       Time snapshot  & $s = 0$  & $s = 0.5$ & $s = 0.9$  \\ 
       \noalign{\smallskip}\hline\noalign{\smallskip}
        ${\rm AE}$ (mean $\pm$ std) &  $0.31 \pm 0.17$ & $0.47 \pm 0.44$ & $0.17 \pm 0.18$ \\ 
        \noalign{\smallskip}
       ${\rm RE}$ (mean $\pm$ std) &  $0.02 \pm 0.01$ & $0.06 \pm 0.06$ & $0.15 \pm 0.17$ \\ 
    \noalign{\smallskip}\hline
    \end{tabular}
    \caption{Average absolute and relative error between $\Phi_{\rm NN}^{(1)}, \ldots, \Phi_{\rm NN}^{(10)}$ and $\Phi_{\rm FEM}$ across all sampled points at different time steps.}
    \label{tab:error_tab }
\end{table}

In~\cref{tab:Phi_J}, we compare the value function approximation for one of the trained models, $\Phi_{\rm NN}^{(1)}$, to the value of the control objective, $J(-\nabla \Phi_{\rm NN})$, at the initial state $\bfx = (-1.5,-1.5)^\top$ and time $t=0$. Since the system dynamics are stochastic, we generate $12,000$ trajectories starting from $\bfx$ using the learned feedback control to calculate the control objective $J$ for each trajectory. We then use the sample average as a proxy for the expected value.
We use a finer step size of $\du s = 0.005$ than the one used in training to get an accurate approximation. We observe that the discrepancy between the value estimate and the actual cost is almost negligible for the FEM solution.
For the neural network approximation, the value estimate is about  4\% smaller than the actual control objective, which indicates that the value estimates can be overly optimistic.

\begin{table}[t]
    \centering
    \begin{tabular}{ccccc}
    \hline  \noalign{\smallskip}
       Initial state  & $\Phi_{\rm FEM}(0)$  & $\Phi_{\rm NN}(0)$ & $J_{\rm FEM}$ & $J_{\rm NN}$  \\ 
       \noalign{\smallskip}\hline\noalign{\smallskip}
        $\bfx_{\text{init}} = (-1.5, -1.5)^\top$ &  $14.67$ & $14.48$ & $14.68$  & $15.33$\\ 
    \noalign{\smallskip}\hline
    \end{tabular}
    \caption{Discrepancy between the value function $\Phi$ and the control objective $J$ at some initial state.}
    \label{tab:Phi_J}
\end{table}

On the hardware used for our experiments, both approaches showed comparable time-to-solution. The neural network training took approximately 20 minutes using the GPU, while the FEM solution was obtained in roughly one hour using the CPU. However, the FEM approach requires a computational mesh, making it infeasible for $d>4$, which is the primary use case for our proposed method.

%%%%%%%%%%%%%%%%%%%%%%%%%%%%%%%%%%%%%%%%%%%%%%%%%%%%%%%%%
\subsection{100-dimensional example}
\label{s:100d-exp}
%%%%%%%%%%%%%%%%%%%%%%%%%%%%%%%%%%%%%%%%%%%%%%%%%%%%%%%%%
We consider the 100-dimensional benchmark SOC problem also used in \cite{Han-2017,Han_BSDE} with  initial state $\bfx = (0,0,\dots,0)^\top\in \R^{100}$ corresponding to time $t=0$.
The drift and diffusion of the system are given by 
\[ f(s,\bfz,\bfu) = 2\bfu\quad \text{and} \quad \sigma = \sqrt{2}, \]
respectively. The terminal and Lagrangian cost  are
\begin{align}\label{eq:100D_term}
	G(\bfz) = \ln \left(\frac{1 + \|\bfz \|^2 }{2}\right),\quad \text{and}\quad L(s,\bfz,\bfu ) =   \|\bfu\|^2,
\end{align} 
respectively.
We compute the Hamiltonian \cref{eq:H_stoc} as
\begin{align*}
	\begin{split}
		H(s,\bfz,\bfp,\bfM) & =\sup_{\bfu\in U}\left\{\frac{\sigma}{2}\text{tr}\left(\bfM\right) + \bfp \cdot f(s,\bfz,\bfu)-L(s,\bfz,\bfu) \right\}\\
		& = \sup_{\bfu\in U}\left\{\frac{1}{\sqrt{2}}\text{tr}\left(\bfM\right) + \bfp \cdot 2\bfu - \|\bfu\|^2 \right\}
	\end{split}
\end{align*}
Using the first-order necessary condition we get
\begin{align*}
	0 = 2\bfu - 2\bfp\implies \bfu = \bfp,
\end{align*}
and using this closed form for $\bfu$, the Hamiltonian is given by
\[ H(s,\bfz,\bfp,\bfM) =\frac{1}{\sqrt{2}} \text{tr}\left(\bfM\right) + \|\bfp\|^2 .\]
Hence, the HJB equation satisfied by the value function, $\Phi(\cdot,\cdot)$, reads
\begin{equation}\label{eq:HJB_100D}
	\frac{\partial }{\partial s} \Phi(s,\bfz) + \Delta \Phi(s,\bfz) - \| \nabla \Phi(s,\bfz) \|^2 = 0, \quad \Phi(T,\bfz) = G(\bfz).
\end{equation}
and its solution is given by
\begin{equation}\label{eq:HJB_100D_sol}
	\Phi(s,\bfz) = -\ln\left(\bbE \left( \exp\left(-G\left(\bfz + \sqrt{2}\,\du W(T-s)\right)\right) \right)\right),
\end{equation}
which we use to test the performance of our method.

Finally, we note that the forward SDE~\cref{eq:sampling} we propose to use for sampling the state space simplifies to
\begin{equation}\label{eq:FBSDE100d_ours}
		\bfz_{i+1} = \bfz_{i}- 2\, \nabla_{\bfz}\Phi(s_i,\bfz_i)\du s + \sqrt{2} \du \bfW_{i}.
\end{equation}

\subsubsection{The importance of sampling}
To demonstrate the impact of using the feedback form to sample the state space, we use the same neural network model as in \cite{raissi2018forward}, which is given by a five-layer feed-forward neural network with 256 neurons per hidden layer to approximate the solution $\Phi(s,\bfz)$. We partition the time interval $[0,1]$ using 50 uniformly spaced points.
We use the same penalty parameters as in the original code, that is, $\beta=(1,0,20,1,1)$.
We use the Adam optimizer  \cite{kingma2014adam} to update the parameters of the network with a batch size of 64 using 50,000 iterations. The average cost per 100 iterations was  $27$s using the CPU. For the following experiments, we use \cref{eq:full_opt} excluding $P_{\rm HJB}$ penalty and compare our method with FBSNNs in \cite{raissi2018forward}.

\begin{figure}[!ht]
        \centering
        \resizebox{\textwidth}{!}{
            \input{figures/100D_problem/comb_results/HJBall.tex}
        }
    \caption{\label{fig:100D_sol}Solution to \cref{eq:HJB_100D} obtained using our method (left column) and the method in \cite{raissi2018forward} (right column)}
\end{figure}

\begin{table}[ht]
    \centering
    \begin{tabular}{ccc}
    \hline  \noalign{\smallskip}
       Method  & \begin{tabular}{cc} \multicolumn{2}{c}{20k iterations}\\
                           RE & RE$_0$
                         \end{tabular} &\begin{tabular}{ll}
                    \multicolumn{2}{c}{50k iterations}\\
                           RE & RE$_0$
                         \end{tabular}  \\ 
       \noalign{\smallskip}\hline\noalign{\smallskip}
        FBSNN &  \begin{tabular}{cc}
                           0.54\% & 0.12\%
                         \end{tabular} 
                         &  \begin{tabular}{cc}
                           0.39\% & 0.045\%
                         \end{tabular} \\ 
        \noalign{\smallskip}
        Ours &  \begin{tabular}{cc}
                           0.48\% & 0.0083\%
                         \end{tabular} 
                         &  \begin{tabular}{cc}
                           0.39\% & 0.012\%
                         \end{tabular} \\ 
    \noalign{\smallskip}\hline
    \end{tabular}
    \caption{Relative errors for \cref{eq:HJB_100D}  obtained using our method and method in \cite{raissi2018forward}}
    \label{tab:rel_HJB }
\end{table}

In \Cref{fig:100D_sol} we plot the exact solution (black-dashed line) \cref{eq:HJB_100D_sol}, the learned solution using our approach (blue-solid line) and the solution learned using FBSNNs (red-solid line) along five random trajectories. In the top row, we present the results obtained after training the networks for 20,000 iterations with a learning rate of $10^{-3}$ and the bottom row presents the results after training the networks for 20K and 30K iterations with learning rates $10^{-3}$ and $10^{-4}$, respectively. These results suggest that our approach approximates the value function better, especially in early iterations, as compared to the FBSNNs. 

In \cref{tab:rel_HJB }, we also compare the learned solutions to the exact solution $\Phi$ in \cref{eq:HJB_100D_sol} by computing the average relative errors,
\[ RE = \frac{\|\Phi(\cdot,\cdot;\bfth) - \Phi(\cdot,\cdot) \|_{2}}{\|\Phi\|_2},\quad   RE_0 = \frac{|\Phi(0,\bfz(0);\bfth) - \Phi(0,\bfz(0)) |}{|\Phi(0,\bfz(0))|},\] 
for ten random trajectories. 
Our method attains lower errors, especially for the initial values and at the earlier iterations.
\subsubsection{Initial states from a distribution}
\begin{figure}
    \centering
    \includegraphics[width=0.5\textwidth]{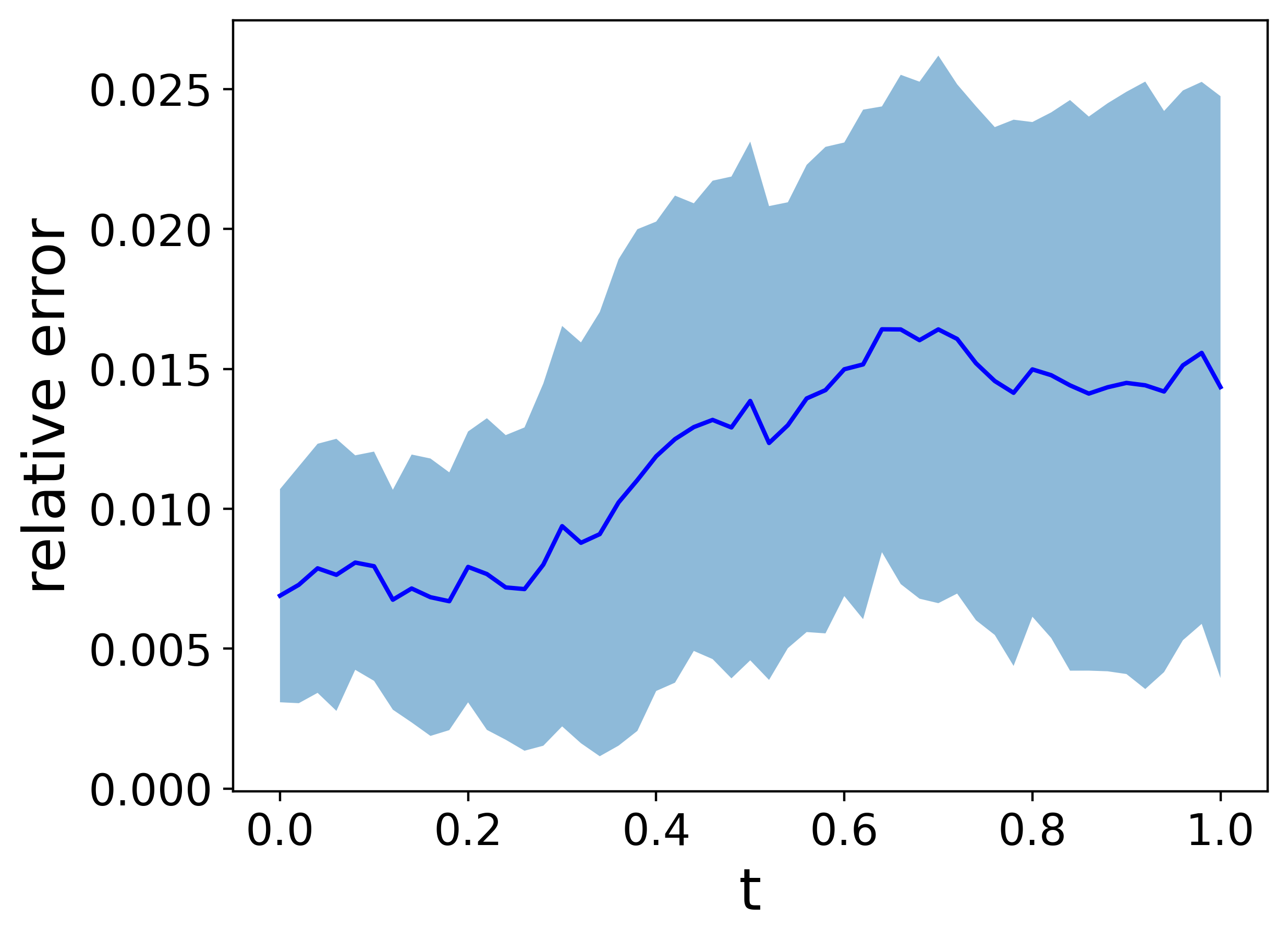}
    \caption{Mean and variance of the errors relative to \cref{eq:HJB_100D_sol} in the learned value function for ten random trajectories obtained by sampling initial states from a distribution using our method after 20k iteration.}
    \label{fig:100d_rho_test}
\end{figure}

We demonstrate the versatility of our method beyond fixed initial states, especially in addressing input states following a given distribution. Specifically, we sample $\bfx$ from a distribution $\rho = \mathcal{N}(\mathbf{0}, 0.5 \cdot \bfI_{100})$.
We repeat the training process with 20k iterations, maintaining the same hyperparameters, but increasing the batch size to $512$ from $64$ and choosing $\beta_3=50$. In \cref{fig:100d_rho_test}, we present the mean and variance of the relative errors of the errors relative to \cref{eq:HJB_100D_sol} in the learned value function for ten random trajectories. As expected to the higher complexity of the problem, the maximum relative error over the time interval increased to $1.5\%$, which is slightly larger than in the original problem.

\subsubsection{Shifted target}\label{s:100d-exp-mod}
In the example above, the minimizer of the terminal function coincides with the initial state $\bfx=(0,0,\dots,0)^\top$. 
Therefore, even a random walk without drift (as used in  \cite{raissi2018forward,Han_BSDE}) will sample around the optimal terminal state, which is critical to accurately approximate the value function. 
This also means that after training using our approach, the drift term in the sampler is relatively small and that the above experiment does not fully show the advantages of our method.

To shed more light on the importance of sampling, we modify the terminal cost to  
$$G(z) = 1000\ln \left(\frac{1 + \|\bfz - \bfz_{\rm target}\|^2 }{2}\right),$$ 
with $\bfz_{\rm target} = (3,3,\dots, 3)^{T}$, so that the target for the state variable $\bfz$ at final time $T$ no longer coincides with the initial state. 
Similar to the two-dimensional test problem in \cref{2D problem}, solving the modified problem now requires sampling around the target and we expect to benefit from the added drift term.

\begin{figure}[t]
	\centering
	\includegraphics[width = 0.5\linewidth]{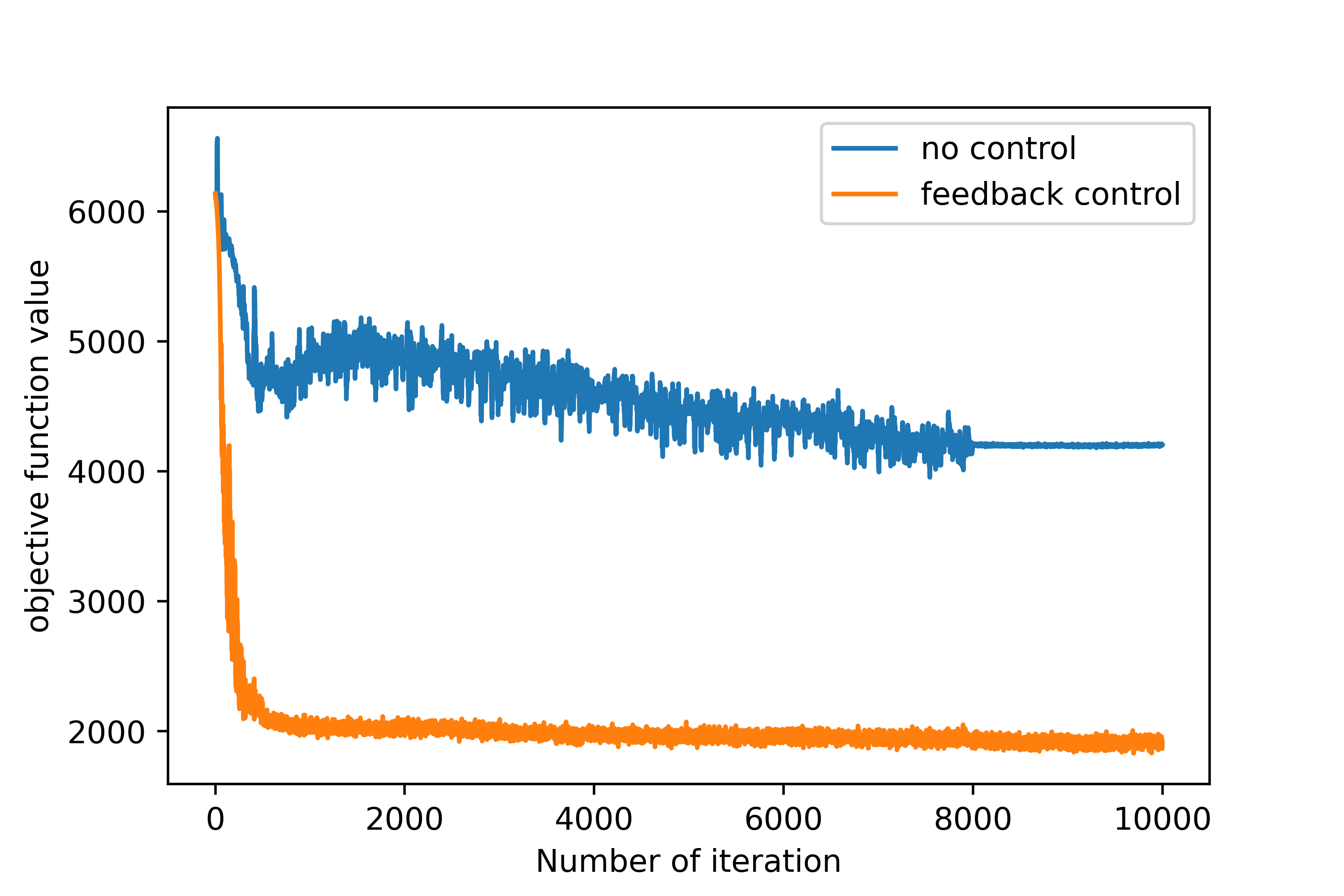}%
	\resizebox{0.36\linewidth}{!}{
            \input{figures/100D_problem/shifted_example.tex}
        }%
	\caption{Computational results for the modified 100-dimensional benchmark problem in~\cref{s:100d-exp-mod}. Left: Control objective for both methods given the same initial state, the blue line represents results using FBSNNs in \cite{raissi2018forward}, and the orange line denotes our method.  Right: Trajectory examples generated using learned value functions on two randomly selected dimensions. The orange line represents our method and blue line FBSNNs.}
	\label{fig:100d_loss}
\end{figure}

We compare our method to FBSNNs on the modified problem keeping the same network structure and hyper-parameters. We use a smaller $\sigma = \frac{2\sqrt{2}}{5}$ to improve training speed. We evaluate the performance of the methods using the objective functional $J$ defined in ~\cref{eq:Joc_stoc} at the control obtained from the feedback form via the respective value function approximations. For this experiment, we use a GPU to train and the results of this comparison are shown in \cref{fig:100d_loss}.

To reduce the effect of the Brownian motion, we run the experiments for each method on the same problem five times and plot the average values corresponding to training iterations. Furthermore, since the primary goal for this example is to explore the difference between sampling strategies, we select much higher weights for the control objective such that we have faster initial convergence for the control variable. 

As can be seen in \cref{fig:100d_loss} (left), our method not only yields faster initial convergence but also achieves a considerably lower control objective. This indicates that the controls obtained from our approach are  more effective, that is, they are closer to optimal. It is also worth pointing out that due to the high terminal cost we assigned when designing the problem, it takes very few iterations to locate the correct state-time region that the optimal solution resides in.
Since FBSNNs use a Brownian motion with no drift, the sampling is unlikely to discover the target.
Consequently, the generated trajectories in \cref{fig:100d_loss} (right)  from our method approximately reach the target, while the trajectories obtained from the FBSNN method stay closer to the initial state. Do note additional hyperparameter tuning and training will be needed if one aims to solve the underlying HJB equation accurately as well.

\subsection{Quadcopter Problem with Nonlinear Dynamics} \label{sec:quadcop}
We test our proposed method's ability to deal with nonlinear dynamics using the stochastic version of the quadcopter trajectory planning problem also considered in \cite{lin2018splitting,onken2020neural}. We sample initial states from a Gaussian distribution centered at $\bfx = [-1.5,-1.5,-1.5,0,\dots,0]^\top$ and set $\bfx_{\text{target}} = [2,2,2,0,\dots,0]^\top$. Here $d=12$, given the state variable $\bfz = [z_1,z_2,z_3,z_4,z_5,z_6,z_7,z_8,z_9,z_{10},z_{11},z_{12}]^\top$ the dynamics read
\begin{equation*}
     f(s,\bfz,\bfu) = \left\{
     \begin{aligned}
       &z_7
       \\
       &z_8
       \\
       &z_9
       \\
       &z_{10}
       \\
       &z_{11}
       \\
       &z_{12}
       \\
       &\frac{u_1}{m} f_7(z_4, z_5, z_6) = \frac{u_1}{m}( \sin(z_4) \sin(z_6) + \cos(z_4) \sin(z_5) \cos(z_6))
       \\
       &\frac{u_1}{m} f_8(z_4, z_5, z_6) = \frac{u_1}{m}( -\cos(z_4) \sin(z_6) + \sin(z_4) \sin(z_5) \cos(z_6))
       \\
       &\frac{u_1}{m} f_9(z_5, z_6) - g = \frac{u_1}{m}(\cos(z_5)\cos(z_6)) - g
       \\
       &u_2
       \\
       &u_3
       \\
       &u_4
       \end{aligned}
     \right. 
 \end{equation*}
The controls for the problem are $\bfu = [u_1,u_2,u_3,u_4]^\top \in \mathbb{R}^{4}$. We assume that both the mass $m$ and gravity $g$ are given. The control objective encompasses $L(\bfu(s,\bfz)) = 2 + \| \bfu(s,\bfz) \|^2$, and $G(\bfz(T)) = 2500 \cdot \| \bfz(T) - \bfx_{\rm target} \|^2$. The feedback form with respect to $\Phi$ for this problem takes the form:
 \begin{equation*}
     \label{eq:feedback_quad}
     \begin{split}
         u_1 = \frac{-1}{2m}\left(f_7 \frac{\partial \Phi}{\partial z_7} + f_8 \frac{\partial \Phi}{\partial z_8} + f_9 \frac{\partial \Phi}{\partial z_9} \right), \;\; \\ 
         u_2 = -\frac{1}{2}\frac{\partial \Phi}{\partial z_{10}}, \;\;
         u_3 = -\frac{1}{2}\frac{\partial \Phi}{\partial z_{11}}, \;\;
         u_4 = -\frac{1}{2}\frac{\partial \Phi}{\partial z_{12}}.
     \end{split}
\end{equation*}
The HJB equation and BSDE can be derived using the feedback form accordingly under \cref{s:prelim}. We choose $\sigma=0.2$ for the problem.

We use the network in \cref{eq:NNArchitecture} featuring two layers and 128 neurons per layer for the ResNet. The penalty term is $\beta=(0.1, 0.1, 1.0, 0.1, 0.1)$. We train the network using 6000 iterations of Adam with a batch size of 128. The learning rate initiates at 0.01 and is halved every 1600 iterations.
Since the dynamics in this example is more complex, we discretize the SDE with 100 equidistant steps between $t=0$ and $T=1$ to enhance accuracy. On average, every training iteration took around 2 seconds on the GPU.

\begin{figure}[t]
    \centering
    \includegraphics[width = 0.55\textwidth]{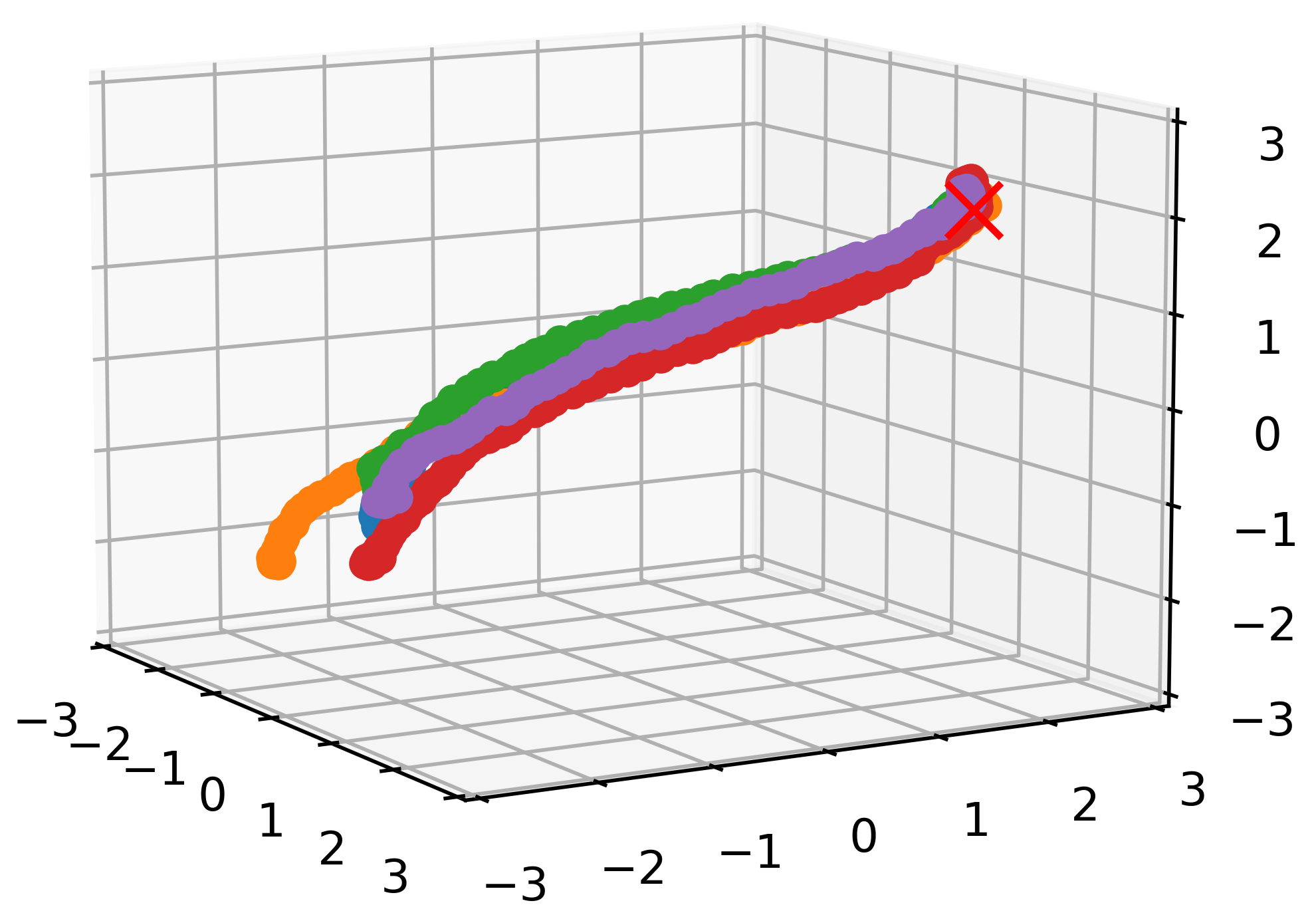}
    \caption{Flight path examples using the learned controller. The target is depicted by a red cross.}
    \label{fig:quad_compile}
\end{figure}

\begin{table}
    \centering
    \begin{tabular}[b]{ccc}\hline
       $J$ at $\bfx$ & evaluated at $\sigma=0.2$ & evaluated at $\sigma=0$\\ \hline
      Deterministic Model & $9.33 \times 10^3$ & $2.18 \times 10^3$\\
      Our Model & $3.34 \times 10^3$  &  -\\
      Accurate $J$ (deterministic) & - & $2.18\times 10^3$ \\ \hline
    \end{tabular}
    \caption{Approximated control objective $J$ for initial state $\bfx = [-1.5,-1.5,-1.5,0,\dots, 0]^\top$. Note the deterministic solution is trained with $\sigma=0$ while ours with $\sigma=0.2$. Value for the accurate solution comes from \cite{onken2020neural}.}
     \label{tab:quad_compile}
\end{table}

The visualization of the trained policy in \cref{fig:quad_compile} shows that the flight trajectories reach the given target from various starting points. 
While we are not aware of an analytical solution to the problem, we compare the performance of our policy to the pre-trained policy from \cite{onken2020neural}, which is trained for the deterministic problem instance.
For each policy, we compute the average value of the control objective over 15,000 randomly chosen trajectories, each using 200 time steps and report the results in \cref{tab:quad_compile}.  As to be expected, while the deterministic solution works well for $\sigma=0$, its performance drops notably when the objective is evaluated with $\sigma=0.2$. 
Since the stochasticity of the dynamics is taken into account during training, our model performs better in this case.

\section{Discussion}
We propose a neural network approach for approximately solving Hamilton-Jacobi-Bellman PDEs arising in high-dimensional stochastic optimal control. 
Similar to existing approaches~\cite{raissi2018forward, Han_BSDE}, we parameterize the value function with a neural network and experiment with different losses to train the network weights.
What sets our work apart from these works is the use of feedback form given by the stochastic Pontryagin maximum principle to design the forward SDE used to explore the state space during training.

Using an intuitive two-dimensional test problem, we visualize that the improved sampling strategy allows us to effectively learn the value function and determine the relevant regions of the state space; see~\cref{2D problem}. 
Based on this insight, we modify the  $100$-dimensional test problem also used in~\cite{raissi2018forward, Han_BSDE} by shifting the minimizer of the terminal costs; see~\cref{s:100d-exp-mod}.
Thereby, we demonstrate that our proposed method dramatically improves the quality of the obtained control.  Using a 12-dimensional quadcopter example whose dynamic is nonlinear in the states, we also demonstrate that our model can handle complicated dynamics; see \cref{sec:quadcop}.

Hyperparameter tuning is crucial in training any neural network model. In our case, choosing different weights for different penalties in the loss function can result in varying outcomes. 
In our experiments, we found that including the control objective, $J$, in the training loss is crucial to obtaining accurate results.
Without this term, we encountered examples where the optimal paths can not be recovered correctly despite having a lower $P_{\rm HJB}$ loss.
One possible cause is that the sampling trajectories, which are also used to compute $P_{\rm HJB}$, do not effectively identify the relevant regions of the state space.
One heuristic we found effective was to use a relatively large weight for the control objective $J$, especially at the beginning of the training. Once the control objective is sufficiently small, we suggest experimenting with the weights corresponding to $P_{\rm BSDE}$ and $P_{\rm HJB}$.

Our numerical experiments also show that the modified forward SDE and the control objective can lead to faster initial convergence compared to the approaches in~\cite{raissi2018forward,Han_BSDE} (see~\cref{s:100d-exp}).

We refer to our method as a semi-global method for solving the HJB equation since we do not aim at estimating the value function accurately globally. 
Instead, we seek to approximate the value function well for states likely to be visited by optimal trajectories of the SOC problem. 
In theory, any point in the state space has a positive probability of being visited due to the stochasticity in the dynamics. 
However, the histogram plots in~\cref{fig:sampling} suggest that the density of the optimal trajectory is concentrated in a small subset of the state space.
Therefore, focusing the exploration on these subsets may have practical advantages over FSDEs that are pure random walks.

Another benefit of our proposed forward SDE compared to purely random exploration is that it coincides with the characteristic curves of the HJB equation as the stochasticity of the system is reduced.
Therefore, our work can be seen as an extension of the neural network approaches for deterministic control problems in~\cite{onken2020neural}.

Compared to neural network approaches for semi-linear elliptic/parabolic PDEs such as~\cite{raissi2018forward, Han_BSDE} it is important to highlight that our approach is limited to  HJB equations arising in stochastic optimal control.
Since our forward SDE is derived from optimality principles, extending it to other high-dimensional PDEs  (for example, Black Scholes and Allen Cahan equations) is not obvious and may be impossible.

In future work, we will apply our approach to problems with non-constant and non-scalar diffusion coefficients.
Although the theoretical framework supports this, we are not aware of any practical algorithms for this case.
Since some SOC problems lead to non-smooth value functions and may include noise terms that depend on the control, testing our scheme on these problems is also a possible extension of our work. 
To overcome the challenge of efficiently computing Hessians of the value function needed for more general HJB penalties, one can use the hessQuik package~\cite{Newman2022}.

\bibliographystyle{siamplain}
\bibliography{refs}
\end{document}

%% file: ex_shared.tex
\usepackage{lipsum}
\usepackage{amsfonts}
\usepackage{graphicx}
\usepackage{epstopdf}
\usepackage{algorithmic}
\ifpdf
  \DeclareGraphicsExtensions{.eps,.pdf,.png,.jpg}
\else
  \DeclareGraphicsExtensions{.eps}
\fi

\newsiamremark{remark}{Remark}
\newsiamremark{hypothesis}{Hypothesis}
\crefname{hypothesis}{Hypothesis}{Hypotheses}
\newsiamthm{claim}{Claim}

\headers{A neural network approach for stochastic optimal control}{X. Li, D. Verma and L. Ruthotto}

\title{A neural network approach for stochastic optimal control\thanks{Submitted to the editors March 10, 2023.
\funding{This work was supported in part by NSF awards DMS 1751636 and DMS 2038118s, AFOSR grant FA9550-20-1-0372, and US DOE Office of Advanced Scientific Computing Research
Field Work Proposal 20-023231.}}}

\author{Xingjian Li\thanks{Department of Mathematics, Emory University, Atlanta, GA
  (\email{xingjian.li@emory.edu}, \\
  \email{deepanshu.verma@emory.edu}, \email{lruthotto@emory.edu}).}
\and Deepanshu Verma\footnotemark[2]
\and Lars Ruthotto\footnotemark[2]}

\usepackage{amsopn}

\usepackage{comment}
\usepackage{subcaption} % for subfigures

\usepackage{latexsym,dsfont}
\usepackage{amsmath,amssymb,amsfonts}
\usepackage[utf8]{inputenc}

\usepackage{xcolor}

\newsiamremark{assumption}{Assumption}

\usepackage{pgfplots}
\usepgfplotslibrary{groupplots}
\pgfplotsset{compat=1.11}
\usepackage{tikz}

%% file: abbrv.tex
\DeclareMathOperator*{\argmax}{arg\,max}

\newcommand{\bfa}{\boldsymbol{a}}
\newcommand{\bfb}{\boldsymbol{b}}

\newcommand{\bfn}{\boldsymbol{n}}

\newcommand{\bfp}{\boldsymbol{p}}

\newcommand{\bfu}{\boldsymbol{u}}

\newcommand{\bfw}{\boldsymbol{w}}
\newcommand{\bfx}{\boldsymbol{x}}
\newcommand{\bfy}{\boldsymbol{y}}
\newcommand{\bfz}{\boldsymbol{z}}

\newcommand{\bfA}{\boldsymbol{A}}

\newcommand{\bfI}{\boldsymbol{I}}

\newcommand{\bfK}{\boldsymbol{K}}

\newcommand{\bfM}{\boldsymbol{M}}

\newcommand{\bfU}{\boldsymbol{U}}
 
\newcommand{\bfW}{\boldsymbol{W}}

\newcommand{\bfZ}{\boldsymbol{Z}}

\newcommand{\bbP}{\mathbb{P}}
\newcommand{\bbF}{\mathbb{F}}

\newcommand{\Fc}{{\mathcal F}}
\newcommand{\Hc}{{\mathcal H}}

\newcommand{\bftheta}{{\boldsymbol \theta}}

\newcommand{\R}{\ensuremath{\mathds{R}}}
\newcommand{\E}{\ensuremath{\mathds{E}}}
\newcommand{\bbE}{\E} % so we all use the same

\def\du{\ensuremath{\mathrm{d}}}

\newcommand{\bfth}{\boldsymbol{\theta}}

%% file: figures/100D_problem/comb_results/HJBall.tex
% This file was created with tikzplotlib v0.10.1.
\begin{tikzpicture}

\definecolor{darkgray176}{RGB}{176,176,176}

\begin{groupplot}[group style={group size=2 by 2,vertical sep=2.3cm}]
\nextgroupplot[
tick align=outside,
tick pos=left,
x grid style={darkgray176},
xmin=-0.05, xmax=1.05,
xtick style={color=black},
y grid style={darkgray176},
ymin=4.31354012447748, ymax=4.87874818304316,
ytick style={color=black},
legend to name=zelda ,legend style={legend columns=3}
]
\addplot [very thick, blue]
table {%
0 4.59004211425781
0.02 4.59153175354004
0.04 4.59182643890381
0.06 4.58934307098389
0.08 4.59127140045166
0.1 4.58570671081543
0.12 4.58411598205566
0.14 4.58827209472656
0.16 4.59220600128174
0.18 4.60179424285889
0.2 4.601731300354
0.22 4.60605907440186
0.24 4.59956026077271
0.26 4.61525011062622
0.28 4.6157922744751
0.3 4.61138916015625
0.32 4.62578105926514
0.34 4.63942718505859
0.36 4.630211353302
0.38 4.61196994781494
0.4 4.62558078765869
0.42 4.62359046936035
0.44 4.62130737304688
0.46 4.60811519622803
0.48 4.61933660507202
0.5 4.62353420257568
0.52 4.62537860870361
0.54 4.63781499862671
0.56 4.61366748809814
0.58 4.60750007629395
0.6 4.59767627716064
0.62 4.61678409576416
0.64 4.63412618637085
0.66 4.62354946136475
0.68 4.62224531173706
0.7 4.61921787261963
0.72 4.58555746078491
0.74 4.55678462982178
0.76 4.54012107849121
0.78 4.54891443252563
0.8 4.56525754928589
0.82 4.59029388427734
0.84 4.56444358825684
0.86 4.58250427246094
0.88 4.58716201782227
0.9 4.60422801971436
0.92 4.61022233963013
0.940000000000001 4.59535598754883
0.960000000000001 4.58325338363647
0.980000000000001 4.56237411499023
1 4.54009389877319
};
\addlegendentry{Learned solution }
\addplot [semithick, black, dashed]
table {%
0 4.58966101118138
0.02 4.5929439424684
0.04 4.59087946362921
0.06 4.58449181426891
0.08 4.58334692657954
0.1 4.57456752133998
0.12 4.5690750132966
0.14 4.57066317628971
0.16 4.57458815307669
0.18 4.59632335406394
0.2 4.5960825037045
0.22 4.59931829860685
0.24 4.59244315774529
0.26 4.61307531656651
0.28 4.61672297028313
0.3 4.6092331841463
0.32 4.63617269422198
0.34 4.64588303070427
0.36 4.62940570060167
0.38 4.61104851169188
0.4 4.62674903471665
0.42 4.60663549863161
0.44 4.60421534993979
0.46 4.59314303867415
0.48 4.60014107274389
0.5 4.61415617392663
0.52 4.62096761413443
0.54 4.62948165995829
0.56 4.6025801356059
0.58 4.59078373088121
0.6 4.59004656923736
0.62 4.62179259811781
0.64 4.64182054091205
0.66 4.61975161451422
0.68 4.6107030503414
0.7 4.60198424522166
0.72 4.56702919967616
0.74 4.54767805133522
0.76 4.52456857647046
0.78 4.53985389374852
0.8 4.5617438584504
0.82 4.58205343779652
0.84 4.56448919851979
0.86 4.56930126095063
0.88 4.5692093054869
0.9 4.60138322847004
0.92 4.60637263459534
0.940000000000001 4.58996857702804
0.960000000000001 4.58095269078082
0.980000000000001 4.55583047507417
1 4.52019028028441
};
\addlegendentry{Exact solution }
\addlegendimage{red}
\addlegendentry{FBSNN solution }
\addplot [semithick, blue]
table {%
0 4.59004211425781
0.02 4.58766746520996
0.04 4.58329963684082
0.06 4.58525943756104
0.08 4.58378124237061
0.1 4.58188056945801
0.12 4.57320690155029
0.14 4.57240390777588
0.16 4.57192802429199
0.18 4.56653451919556
0.2 4.56520843505859
0.22 4.57002258300781
0.24 4.57322835922241
0.26 4.56620454788208
0.28 4.56946659088135
0.3 4.56687831878662
0.32 4.57180213928223
0.34 4.5819730758667
0.36 4.57146167755127
0.38 4.5695686340332
0.4 4.57300233840942
0.42 4.57969665527344
0.44 4.58504962921143
0.46 4.58958435058594
0.48 4.58995008468628
0.5 4.58948516845703
0.52 4.58362865447998
0.54 4.57050466537476
0.56 4.57644367218018
0.58 4.55684995651245
0.6 4.57404661178589
0.62 4.60435199737549
0.64 4.60392665863037
0.66 4.62096214294434
0.68 4.60294628143311
0.7 4.62680673599243
0.72 4.62099361419678
0.74 4.59514999389648
0.76 4.63342142105103
0.78 4.61055612564087
0.8 4.59315347671509
0.82 4.5354642868042
0.84 4.53445243835449
0.86 4.50140380859375
0.88 4.46319341659546
0.9 4.48433876037598
0.92 4.40655946731567
0.940000000000001 4.4253249168396
0.960000000000001 4.42884111404419
0.980000000000001 4.40323925018311
1 4.3760871887207
};
\addplot [semithick, black, dashed]
table {%
0 4.58966101118138
0.02 4.58830112791595
0.04 4.58549251045589
0.06 4.58544739451419
0.08 4.59027715077564
0.1 4.58460265116869
0.12 4.57898906374502
0.14 4.58732074222064
0.16 4.57526864251209
0.18 4.55529155557247
0.2 4.5586873677927
0.22 4.55515055884199
0.24 4.55891168982144
0.26 4.54193308384756
0.28 4.54987590926652
0.3 4.54376006348841
0.32 4.55233702406607
0.34 4.58380429797238
0.36 4.57213465546358
0.38 4.56948099911046
0.4 4.57118112393203
0.42 4.57703032640414
0.44 4.59116143625768
0.46 4.60842249210101
0.48 4.60408161329696
0.5 4.59405608983139
0.52 4.58812582111078
0.54 4.56194956521266
0.56 4.56667933105761
0.58 4.53668716951068
0.6 4.55743136496276
0.62 4.59130306532069
0.64 4.59264425931531
0.66 4.60532127180718
0.68 4.58325629856536
0.7 4.62150500515393
0.72 4.6283239261581
0.74 4.58387569140499
0.76 4.62697609525391
0.78 4.60090518066718
0.8 4.57625021770942
0.82 4.52907560304646
0.84 4.51596856313617
0.86 4.47988827521801
0.88 4.43636947354865
0.9 4.4660760626034
0.92 4.3890421029981
0.940000000000001 4.41485982901935
0.960000000000001 4.4151917846631
0.980000000000001 4.38697142800921
1 4.33923139986683
};
\addplot [semithick, blue]
table {%
0 4.59004211425781
0.02 4.58716249465942
0.04 4.59125423431396
0.06 4.58919143676758
0.08 4.58484315872192
0.1 4.58379983901978
0.12 4.58279085159302
0.14 4.57856369018555
0.16 4.58891105651855
0.18 4.59556341171265
0.2 4.59228944778442
0.22 4.59540939331055
0.24 4.5985312461853
0.26 4.59141540527344
0.28 4.60421276092529
0.3 4.60481214523315
0.32 4.60890054702759
0.34 4.58737373352051
0.36 4.57932472229004
0.38 4.58122587203979
0.4 4.57918834686279
0.42 4.57075023651123
0.44 4.56157922744751
0.46 4.57640647888184
0.48 4.58907890319824
0.5 4.59944009780884
0.52 4.58444786071777
0.54 4.58341789245605
0.56 4.60724782943726
0.58 4.6051664352417
0.6 4.61759090423584
0.62 4.60643672943115
0.64 4.59785938262939
0.66 4.58434438705444
0.68 4.58900928497314
0.7 4.60171604156494
0.72 4.60014581680298
0.74 4.54998111724854
0.76 4.49330520629883
0.78 4.50129890441895
0.8 4.50158214569092
0.82 4.53604984283447
0.84 4.53465986251831
0.86 4.57369756698608
0.88 4.56773042678833
0.9 4.57119417190552
0.92 4.58416175842285
0.940000000000001 4.60953187942505
0.960000000000001 4.58758115768433
0.980000000000001 4.60445213317871
1 4.66053342819214
};
\addplot [semithick, black, dashed]
table {%
0 4.58966101118138
0.02 4.59028637968382
0.04 4.59592822179054
0.06 4.58954031275187
0.08 4.59010859530629
0.1 4.59187874342429
0.12 4.58850064068677
0.14 4.5787706625774
0.16 4.5895023157562
0.18 4.6136845318094
0.2 4.61633846781889
0.22 4.61809401318038
0.24 4.62849963648102
0.26 4.61076823029989
0.28 4.63573827298917
0.3 4.63525234285093
0.32 4.64800644146972
0.34 4.60962568023272
0.36 4.61021040309526
0.38 4.6254606172968
0.4 4.63026438640842
0.42 4.61209118784059
0.44 4.59241997978554
0.46 4.61351101047098
0.48 4.62557881089843
0.5 4.64344028913608
0.52 4.61556485435589
0.54 4.60964047639441
0.56 4.62741100280267
0.58 4.61632227323997
0.6 4.62635444182949
0.62 4.61163361707848
0.64 4.59485395788143
0.66 4.58099711578488
0.68 4.59810362363056
0.7 4.60696421444159
0.72 4.59693545873318
0.74 4.52680342274356
0.76 4.4820381074773
0.78 4.49901478191169
0.8 4.49856914825485
0.82 4.53185574932943
0.84 4.53613823921214
0.86 4.5575604976296
0.88 4.5448630572785
0.9 4.55483712990372
0.92 4.56405789900923
0.940000000000001 4.58894374464275
0.960000000000001 4.5657032631485
0.980000000000001 4.57767618451131
1 4.63248264251707
};
\addplot [semithick, blue]
table {%
0 4.59004211425781
0.02 4.59880971908569
0.04 4.6007080078125
0.06 4.59620475769043
0.08 4.59208583831787
0.1 4.59035110473633
0.12 4.59152221679688
0.14 4.5816125869751
0.16 4.58731842041016
0.18 4.58240699768066
0.2 4.58564710617065
0.22 4.58977603912354
0.24 4.59241199493408
0.26 4.59264087677002
0.28 4.59736442565918
0.3 4.60105991363525
0.32 4.61855316162109
0.34 4.62293434143066
0.36 4.62837409973145
0.38 4.63296747207642
0.4 4.62808513641357
0.42 4.64242649078369
0.44 4.66589450836182
0.46 4.69818925857544
0.48 4.707688331604
0.5 4.73394393920898
0.52 4.70007562637329
0.54 4.71290731430054
0.56 4.72689914703369
0.58 4.71810817718506
0.6 4.7335991859436
0.62 4.71552753448486
0.64 4.69647312164307
0.66 4.66696214675903
0.68 4.68079280853271
0.7 4.68350315093994
0.72 4.71307182312012
0.74 4.73338794708252
0.76 4.71756649017334
0.78 4.69082880020142
0.8 4.73747396469116
0.82 4.76893615722656
0.84 4.77522087097168
0.86 4.81543159484863
0.88 4.85305690765381
0.9 4.8394889831543
0.92 4.81767749786377
0.940000000000001 4.78114032745361
0.960000000000001 4.77295017242432
0.980000000000001 4.78579092025757
1 4.75862789154053
};
\addplot [semithick, black, dashed]
table {%
0 4.58966101118138
0.02 4.59157324741338
0.04 4.59688708002018
0.06 4.59811002387039
0.08 4.60627393552665
0.1 4.60120755445453
0.12 4.60526617095578
0.14 4.58128190442637
0.16 4.59639172744734
0.18 4.58743877296565
0.2 4.59250523797447
0.22 4.60052171676879
0.24 4.6032637364095
0.26 4.60428816234575
0.28 4.60496211886787
0.3 4.62010469035358
0.32 4.64095480746894
0.34 4.64263908330756
0.36 4.64288328849799
0.38 4.64402324332057
0.4 4.62881026208035
0.42 4.64290685718239
0.44 4.67780768006253
0.46 4.7234237408651
0.48 4.73302560559588
0.5 4.76373035118575
0.52 4.71516728464007
0.54 4.72908692664465
0.56 4.7457765435331
0.58 4.7343843997615
0.6 4.75011132449809
0.62 4.72063938687201
0.64 4.70430155891838
0.66 4.66888718543811
0.68 4.68617950483006
0.7 4.6935472057591
0.72 4.71474249269246
0.74 4.73065224296676
0.76 4.72868073254464
0.78 4.69699835309954
0.8 4.7339272606907
0.82 4.76025846607367
0.84 4.74646076050267
0.86 4.79354065838749
0.88 4.83193955189723
0.9 4.82137363023095
0.92 4.79303008924636
0.940000000000001 4.75333493407495
0.960000000000001 4.72278835661518
0.980000000000001 4.71514199015369
1 4.67876550244397
};
\addplot [semithick, blue]
table {%
0 4.59004211425781
0.02 4.59014081954956
0.04 4.58951473236084
0.06 4.59260129928589
0.08 4.58907413482666
0.1 4.59095907211304
0.12 4.59193992614746
0.14 4.593994140625
0.16 4.59164047241211
0.18 4.60370826721191
0.2 4.6041898727417
0.22 4.60864067077637
0.24 4.59375047683716
0.26 4.590989112854
0.28 4.57461166381836
0.3 4.58894109725952
0.32 4.57624816894531
0.34 4.56987190246582
0.36 4.56753158569336
0.38 4.56766223907471
0.4 4.5593695640564
0.42 4.57636260986328
0.44 4.55752992630005
0.46 4.54619884490967
0.48 4.53831434249878
0.5 4.543701171875
0.52 4.53470611572266
0.54 4.52000284194946
0.56 4.5369987487793
0.58 4.5009708404541
0.6 4.49512100219727
0.62 4.50506210327148
0.64 4.49531745910645
0.66 4.50687503814697
0.68 4.49172687530518
0.7 4.52836084365845
0.72 4.50534629821777
0.74 4.50406932830811
0.76 4.47842311859131
0.78 4.47790145874023
0.8 4.53042221069336
0.82 4.56824541091919
0.84 4.58175992965698
0.86 4.52051067352295
0.88 4.50758743286133
0.9 4.49364566802979
0.92 4.49744510650635
0.940000000000001 4.51086854934692
0.960000000000001 4.61139011383057
0.980000000000001 4.54844188690186
1 4.51711082458496
};
\addplot [semithick, black, dashed]
table {%
0 4.58966101118138
0.02 4.59265273792626
0.04 4.58474053579815
0.06 4.58128903575765
0.08 4.57745156852044
0.1 4.58089734754179
0.12 4.58257079872554
0.14 4.5913135963412
0.16 4.58440486246066
0.18 4.6114204999316
0.2 4.60284319209387
0.22 4.61441692261279
0.24 4.58388106870877
0.26 4.57642868456189
0.28 4.56113736553954
0.3 4.57555301177632
0.32 4.54457662134853
0.34 4.54388253242042
0.36 4.55320668090534
0.38 4.55366472674835
0.4 4.53074550382724
0.42 4.559135691971
0.44 4.53432704230091
0.46 4.53052783952621
0.48 4.51837858521473
0.5 4.53452878711991
0.52 4.52975950824981
0.54 4.51252611089667
0.56 4.53750923370349
0.58 4.49349336954935
0.6 4.48777138114173
0.62 4.51030297636051
0.64 4.49885077683184
0.66 4.51876418058595
0.68 4.5033212438416
0.7 4.5409776060956
0.72 4.50574538528992
0.74 4.50607102003518
0.76 4.48702992271468
0.78 4.48173004731961
0.8 4.53508489112334
0.82 4.57496943460973
0.84 4.57847838117674
0.86 4.53140783635281
0.88 4.51877495822152
0.9 4.50378641024499
0.92 4.50174009094135
0.940000000000001 4.50898455894129
0.960000000000001 4.59992447969104
0.980000000000001 4.52722452627211
1 4.50318761468228
};

\nextgroupplot[
tick align=outside,
tick pos=left,
x grid style={darkgray176},
xmin=-0.05, xmax=1.05,
xtick style={color=black},
y grid style={darkgray176},
ymin=4.29173229184724, ymax=4.79280549161238,
ytick style={color=black}
]
\addplot [semithick, red]
table {%
0 4.58372068405151
0.02 4.58568477630615
0.04 4.58507442474365
0.06 4.58793115615845
0.08 4.58721733093262
0.1 4.59056425094604
0.12 4.59281063079834
0.14 4.59494686126709
0.16 4.5980920791626
0.18 4.60402345657349
0.2 4.60175848007202
0.22 4.60758781433105
0.24 4.61071491241455
0.26 4.60009908676147
0.28 4.59760999679565
0.3 4.58254957199097
0.32 4.59689044952393
0.34 4.58830451965332
0.36 4.58698654174805
0.38 4.60137510299683
0.4 4.59475898742676
0.42 4.57748889923096
0.44 4.58334875106812
0.46 4.56468343734741
0.48 4.54984188079834
0.5 4.54115581512451
0.52 4.53091907501221
0.54 4.51863765716553
0.56 4.5255389213562
0.58 4.54720497131348
0.6 4.56948280334473
0.62 4.58897018432617
0.64 4.57000923156738
0.66 4.5417652130127
0.68 4.53551626205444
0.7 4.4994740486145
0.72 4.4817476272583
0.74 4.48267650604248
0.76 4.46094989776611
0.78 4.48984384536743
0.8 4.51263666152954
0.82 4.51676750183105
0.84 4.51610374450684
0.86 4.52374458312988
0.88 4.51395893096924
0.9 4.55009412765503
0.92 4.54142713546753
0.940000000000001 4.52361011505127
0.960000000000001 4.54143810272217
0.980000000000001 4.52712678909302
1 4.55265140533447
};

\addplot [semithick, black, dashed]
table {%
0 4.58966101118138
0.02 4.59168554934643
0.04 4.59114714885045
0.06 4.59736878856505
0.08 4.60150819679974
0.1 4.6174602467402
0.12 4.6184761378316
0.14 4.62625442148295
0.16 4.63070186047094
0.18 4.64575428316228
0.2 4.63307813526105
0.22 4.63627501381225
0.24 4.64293142010073
0.26 4.61376272375088
0.28 4.60364949377519
0.3 4.58774561263564
0.32 4.61775237281303
0.34 4.60239964552265
0.36 4.59765648595001
0.38 4.62012192236195
0.4 4.61172885514909
0.42 4.59096251201627
0.44 4.60062084865938
0.46 4.56848627532736
0.48 4.54522573797858
0.5 4.52703204496042
0.52 4.50168990473442
0.54 4.49286494140208
0.56 4.50035856167028
0.58 4.5341744279782
0.6 4.55775132540084
0.62 4.58747568218259
0.64 4.56397130886382
0.66 4.52657523278522
0.68 4.52721037136733
0.7 4.48892692239495
0.72 4.45968617970736
0.74 4.46614813795766
0.76 4.44307721239008
0.78 4.48759167762656
0.8 4.51617590504358
0.82 4.5180497257361
0.84 4.5211511132248
0.86 4.53200073443817
0.88 4.52626891863499
0.9 4.57705072148817
0.92 4.56277725602872
0.940000000000001 4.54395755626738
0.960000000000001 4.56312582752069
0.980000000000001 4.55001312510312
1 4.57196645980006
};
\addplot [semithick, red]
table {%
0 4.58372068405151
0.02 4.58177089691162
0.04 4.58256721496582
0.06 4.58191013336182
0.08 4.58212852478027
0.1 4.57919216156006
0.12 4.58092355728149
0.14 4.57535839080811
0.16 4.57105779647827
0.18 4.56569004058838
0.2 4.57009935379028
0.22 4.56376934051514
0.24 4.56884622573853
0.26 4.56688976287842
0.28 4.56546878814697
0.3 4.57076025009155
0.32 4.56060695648193
0.34 4.56277179718018
0.36 4.57220840454102
0.38 4.55961608886719
0.4 4.56705904006958
0.42 4.56143093109131
0.44 4.55336570739746
0.46 4.55321788787842
0.48 4.55253839492798
0.5 4.54552841186523
0.52 4.56350517272949
0.54 4.57185649871826
0.56 4.56878280639648
0.58 4.58206558227539
0.6 4.59861755371094
0.62 4.58327960968018
0.64 4.61485195159912
0.66 4.61847352981567
0.68 4.60616016387939
0.7 4.60013389587402
0.72 4.64123249053955
0.74 4.64351463317871
0.76 4.63353824615479
0.78 4.67361259460449
0.8 4.68106031417847
0.82 4.68519687652588
0.84 4.63987255096436
0.86 4.64052295684814
0.88 4.6528491973877
0.9 4.65151262283325
0.92 4.66290283203125
0.940000000000001 4.63485097885132
0.960000000000001 4.62460517883301
0.980000000000001 4.66116571426392
1 4.69106388092041
};
\addplot [semithick, black, dashed]
table {%
0 4.58966101118138
0.02 4.58732190336552
0.04 4.58798545100439
0.06 4.59457684604707
0.08 4.59468061459037
0.1 4.59661279041461
0.12 4.59937146253054
0.14 4.58290773317163
0.16 4.57871072327621
0.18 4.57413526700598
0.2 4.58156276955015
0.22 4.57256302120542
0.24 4.5766906082294
0.26 4.58253057844265
0.28 4.57887436795094
0.3 4.59818089631187
0.32 4.58375896742432
0.34 4.57781425993581
0.36 4.59967726033333
0.38 4.58341467004901
0.4 4.59894291105675
0.42 4.59169065145376
0.44 4.58030079853549
0.46 4.58560617283259
0.48 4.58002840328237
0.5 4.57748668655518
0.52 4.60303179097919
0.54 4.61466867507062
0.56 4.60894203640534
0.58 4.63541818041985
0.6 4.65478087295446
0.62 4.63220819737764
0.64 4.67160015496832
0.66 4.67943749183167
0.68 4.65988201286098
0.7 4.65057264514072
0.72 4.69935282767488
0.74 4.71133955223958
0.76 4.71303050387577
0.78 4.74528572065837
0.8 4.76935086720186
0.82 4.7700294370776
0.84 4.71538905569149
0.86 4.70459729701442
0.88 4.71864512248508
0.9 4.72662131914593
0.92 4.74294105210646
0.940000000000001 4.70870232968849
0.960000000000001 4.69362913135466
0.980000000000001 4.73257916007111
1 4.75700964805685
};
\addplot [semithick, red]
table {%
0 4.58372068405151
0.02 4.58358764648438
0.04 4.58963775634766
0.06 4.59004020690918
0.08 4.59261322021484
0.1 4.59239387512207
0.12 4.59662246704102
0.14 4.59439849853516
0.16 4.59540128707886
0.18 4.60044002532959
0.2 4.59193229675293
0.22 4.58957672119141
0.24 4.58398056030273
0.26 4.59427690505981
0.28 4.58812618255615
0.3 4.59785270690918
0.32 4.59391117095947
0.34 4.59322738647461
0.36 4.59136295318604
0.38 4.59002876281738
0.4 4.58645582199097
0.42 4.58878898620605
0.44 4.59485244750977
0.46 4.59741878509521
0.48 4.59210824966431
0.5 4.58423948287964
0.52 4.58973693847656
0.54 4.58027076721191
0.56 4.55808067321777
0.58 4.56075859069824
0.6 4.56063413619995
0.62 4.55262994766235
0.64 4.54972696304321
0.66 4.53320455551147
0.68 4.48298072814941
0.7 4.516441822052
0.72 4.5195517539978
0.74 4.53517436981201
0.76 4.53983545303345
0.78 4.52940464019775
0.8 4.50970363616943
0.82 4.48925590515137
0.84 4.45796918869019
0.86 4.4564790725708
0.88 4.42170715332031
0.9 4.38195133209229
0.92 4.43246221542358
0.940000000000001 4.39675903320312
0.960000000000001 4.3883638381958
0.980000000000001 4.38791370391846
1 4.36819458007812
};
\addplot [semithick, black, dashed]
table {%
0 4.58966101118138
0.02 4.59217572435319
0.04 4.59586516561154
0.06 4.59192823977113
0.08 4.59697234254217
0.1 4.60227445250023
0.12 4.61096442011321
0.14 4.61214601026326
0.16 4.60797535329083
0.18 4.60661075474493
0.2 4.59490690170423
0.22 4.58698710900683
0.24 4.57909550037571
0.26 4.60135254812314
0.28 4.59113787818662
0.3 4.60321460823478
0.32 4.59647505349109
0.34 4.58672589072227
0.36 4.58141806284406
0.38 4.57370125879449
0.4 4.56600328369507
0.42 4.57555635868058
0.44 4.59407577059658
0.46 4.59644207442235
0.48 4.5872958961661
0.5 4.58081046601734
0.52 4.58811900035823
0.54 4.57317060707338
0.56 4.53920021700149
0.58 4.53563504381043
0.6 4.54102640221954
0.62 4.52797526739673
0.64 4.52405401704066
0.66 4.50384171961588
0.68 4.4341257898412
0.7 4.47436920940257
0.72 4.47828651779791
0.74 4.49648379425169
0.76 4.5040681913912
0.78 4.48779760341185
0.8 4.46743793241238
0.82 4.43677100028675
0.84 4.40872996613058
0.86 4.41237260116627
0.88 4.36717829213186
0.9 4.3323389815476
0.92 4.39144166167334
0.940000000000001 4.35794337831359
0.960000000000001 4.35068551220178
0.980000000000001 4.34067732208561
1 4.31450834638202
};
\addplot [semithick, red]
table {%
0 4.58372068405151
0.02 4.5832462310791
0.04 4.5817699432373
0.06 4.57852458953857
0.08 4.57665681838989
0.1 4.57530498504639
0.12 4.57240915298462
0.14 4.57201814651489
0.16 4.56981515884399
0.18 4.57819652557373
0.2 4.57888650894165
0.22 4.58464956283569
0.24 4.57626724243164
0.26 4.57616138458252
0.28 4.56698560714722
0.3 4.57516860961914
0.32 4.56716537475586
0.34 4.58499002456665
0.36 4.58967971801758
0.38 4.60137796401978
0.4 4.60426330566406
0.42 4.6164493560791
0.44 4.5858268737793
0.46 4.57882785797119
0.48 4.58740234375
0.5 4.57015037536621
0.52 4.56832599639893
0.54 4.58430576324463
0.56 4.58630657196045
0.58 4.56601285934448
0.6 4.59039688110352
0.62 4.56658363342285
0.64 4.5438928604126
0.66 4.55340051651001
0.68 4.51847362518311
0.7 4.52007579803467
0.72 4.49891757965088
0.74 4.51270961761475
0.76 4.52612590789795
0.78 4.54495191574097
0.8 4.52950477600098
0.82 4.56979274749756
0.84 4.6034984588623
0.86 4.590651512146
0.88 4.58232593536377
0.9 4.60327625274658
0.92 4.6080379486084
0.940000000000001 4.6361985206604
0.960000000000001 4.66364479064941
0.980000000000001 4.66345977783203
1 4.66940069198608
};
\addplot [semithick, black, dashed]
table {%
0 4.58966101118138
0.02 4.58807749414542
0.04 4.58554444929648
0.06 4.575928690645
0.08 4.57108775004133
0.1 4.57134472453665
0.12 4.57146721427222
0.14 4.5762304363108
0.16 4.56605913427115
0.18 4.57886977057525
0.2 4.58278465247121
0.22 4.59591558071833
0.24 4.5798332266534
0.26 4.5800082815063
0.28 4.55981940285525
0.3 4.57702243105795
0.32 4.56490117955918
0.34 4.59404205512029
0.36 4.59981984170041
0.38 4.61624398655899
0.4 4.61806432732703
0.42 4.641554244282
0.44 4.58907083267868
0.46 4.57843700914539
0.48 4.58749504521995
0.5 4.5614348384734
0.52 4.56692970044395
0.54 4.59242857112403
0.56 4.59609624829505
0.58 4.5696610035719
0.6 4.60327402742237
0.62 4.57558929695129
0.64 4.54792825062848
0.66 4.5594583409834
0.68 4.52237541348346
0.7 4.51900559737691
0.72 4.49370634433617
0.74 4.50954620567145
0.76 4.52846533748407
0.78 4.55432989785936
0.8 4.54097399588025
0.82 4.58274003795676
0.84 4.61794645482859
0.86 4.60107015282897
0.88 4.59843465181238
0.9 4.62321134321955
0.92 4.62425672460076
0.940000000000001 4.64844176823905
0.960000000000001 4.67771020713876
0.980000000000001 4.6719885949902
1 4.66958231102988
};
\addplot [semithick, red]
table {%
0 4.58372068405151
0.02 4.58430576324463
0.04 4.58566284179688
0.06 4.58491039276123
0.08 4.58312606811523
0.1 4.58024311065674
0.12 4.57338953018188
0.14 4.56637191772461
0.16 4.56049346923828
0.18 4.55726051330566
0.2 4.55214262008667
0.22 4.55751895904541
0.24 4.54601001739502
0.26 4.54877758026123
0.28 4.54295444488525
0.3 4.52885627746582
0.32 4.53363037109375
0.34 4.52077865600586
0.36 4.51370477676392
0.38 4.48759937286377
0.4 4.48356056213379
0.42 4.47034549713135
0.44 4.47825050354004
0.46 4.5037088394165
0.48 4.50653028488159
0.5 4.50822114944458
0.52 4.50613021850586
0.54 4.51064920425415
0.56 4.47996759414673
0.58 4.49469614028931
0.6 4.47538661956787
0.62 4.47062492370605
0.64 4.50502443313599
0.66 4.47548484802246
0.68 4.48887395858765
0.7 4.48341655731201
0.72 4.47367191314697
0.74 4.46806049346924
0.76 4.45130348205566
0.78 4.44155168533325
0.8 4.48296356201172
0.82 4.503258228302
0.84 4.5012092590332
0.86 4.48702478408813
0.88 4.49303913116455
0.9 4.48480224609375
0.92 4.51196193695068
0.940000000000001 4.55521678924561
0.960000000000001 4.64060258865356
0.980000000000001 4.67690849304199
1 4.708655834198
};
\addplot [semithick, black, dashed]
table {%
0 4.58966101118138
0.02 4.58659437056441
0.04 4.58685745595238
0.06 4.58712326525528
0.08 4.58451837578947
0.1 4.57930549496727
0.12 4.55980569594805
0.14 4.54837583029376
0.16 4.54792580074373
0.18 4.54166650337728
0.2 4.5327892904923
0.22 4.55388166382939
0.24 4.52878052688498
0.26 4.53694736931401
0.28 4.52040125496971
0.3 4.48886549015168
0.32 4.50480841952395
0.34 4.48816495901333
0.36 4.48451639316193
0.38 4.44509982819462
0.4 4.44681560778944
0.42 4.42663018684328
0.44 4.44588710001753
0.46 4.48384219346642
0.48 4.49634953681817
0.5 4.49868760179416
0.52 4.49412384185181
0.54 4.50265696485545
0.56 4.45986476581386
0.58 4.47550256715988
0.6 4.4571174593459
0.62 4.45463402730446
0.64 4.4976372140742
0.66 4.46507772290511
0.68 4.48231677930458
0.7 4.47954990198189
0.72 4.47971880386352
0.74 4.47857155551823
0.76 4.46791933550727
0.78 4.44976578625203
0.8 4.48797493540904
0.82 4.50732624951136
0.84 4.51252278219828
0.86 4.49776543451162
0.88 4.50474477238104
0.9 4.50042142862316
0.92 4.51099777975907
0.940000000000001 4.53710950033209
0.960000000000001 4.61463204956937
0.980000000000001 4.6474161785944
1 4.66878235530101
};
\nextgroupplot[
tick align=outside,
tick pos=left,
x grid style={darkgray176},
xmin=-0.05, xmax=1.05,
xtick style={color=black},
y grid style={darkgray176},
ymin=4.20057925430235, ymax=4.87512165513398,
ytick style={color=black}
]
\addplot [semithick, blue]
table {%
0 4.59025430679321
0.02 4.59163522720337
0.04 4.58825635910034
0.06 4.58759069442749
0.08 4.58417940139771
0.1 4.58854341506958
0.12 4.58258247375488
0.14 4.58115673065186
0.16 4.58112621307373
0.18 4.58371114730835
0.2 4.5839581489563
0.22 4.58808088302612
0.24 4.58481073379517
0.26 4.58086729049683
0.28 4.58338499069214
0.3 4.59329319000244
0.32 4.58459997177124
0.34 4.58995962142944
0.36 4.60026741027832
0.38 4.57676839828491
0.4 4.57889223098755
0.42 4.57317209243774
0.44 4.58189249038696
0.46 4.6014575958252
0.48 4.59211683273315
0.5 4.56922101974487
0.52 4.57057666778564
0.54 4.56690454483032
0.56 4.55265760421753
0.58 4.53168535232544
0.6 4.52563142776489
0.62 4.50369596481323
0.64 4.49566030502319
0.66 4.45824241638184
0.68 4.44608449935913
0.7 4.47219896316528
0.72 4.46291065216064
0.74 4.41301488876343
0.76 4.38535022735596
0.78 4.37434053421021
0.8 4.36608839035034
0.82 4.38417291641235
0.84 4.37742805480957
0.86 4.36768341064453
0.88 4.36349058151245
0.9 4.36328935623169
0.92 4.35755157470703
0.940000000000001 4.29784965515137
0.960000000000001 4.29493570327759
0.980000000000001 4.25021028518677
1 4.23124027252197
};
\addplot [semithick, black, dashed]
table {%
0 4.58966101118138
0.02 4.592686932812
0.04 4.58381043588388
0.06 4.58475275578693
0.08 4.57725601759544
0.1 4.59408839015186
0.12 4.57771442477148
0.14 4.57631043923431
0.16 4.57347795315384
0.18 4.58600155100774
0.2 4.59079512278467
0.22 4.60089204998259
0.24 4.5951041503769
0.26 4.58688531574931
0.28 4.59044690466885
0.3 4.6149677286542
0.32 4.59843736953852
0.34 4.60805465937593
0.36 4.62251407089896
0.38 4.5739382331466
0.4 4.57989487115412
0.42 4.5662444406774
0.44 4.58048162777855
0.46 4.61312424888499
0.48 4.59578523429721
0.5 4.56108191001955
0.52 4.56715120396313
0.54 4.5652740500109
0.56 4.54609631208884
0.58 4.51621167534611
0.6 4.50976422389117
0.62 4.48118333484116
0.64 4.47548018722778
0.66 4.42631750099361
0.68 4.4100282500079
0.7 4.4487228474466
0.72 4.440659805512
0.74 4.38151501952938
0.76 4.35140271215432
0.78 4.34385105232099
0.8 4.34083967430321
0.82 4.36855023624517
0.84 4.36299304354191
0.86 4.35608563350213
0.88 4.35460104171128
0.9 4.35724440891058
0.92 4.35526680553813
0.940000000000001 4.29519480140931
0.960000000000001 4.29726977728851
0.980000000000001 4.25450917890389
1 4.23697608199334
};
\addplot [semithick, blue]
table {%
0 4.59025430679321
0.02 4.58943557739258
0.04 4.58917236328125
0.06 4.58943462371826
0.08 4.59115266799927
0.1 4.59004545211792
0.12 4.58509206771851
0.14 4.58660221099854
0.16 4.58595848083496
0.18 4.58047151565552
0.2 4.57788467407227
0.22 4.58330202102661
0.24 4.57455968856812
0.26 4.5792407989502
0.28 4.5791220664978
0.3 4.58504152297974
0.32 4.59232521057129
0.34 4.59457731246948
0.36 4.60285997390747
0.38 4.61067724227905
0.4 4.60705518722534
0.42 4.60429716110229
0.44 4.58866214752197
0.46 4.5969500541687
0.48 4.61742210388184
0.5 4.63310194015503
0.52 4.64065074920654
0.54 4.65258932113647
0.56 4.67809963226318
0.58 4.69946241378784
0.6 4.69989538192749
0.62 4.69341993331909
0.64 4.73423385620117
0.66 4.75549173355103
0.68 4.76510763168335
0.7 4.75545310974121
0.72 4.75538015365601
0.74 4.80239915847778
0.76 4.78395318984985
0.78 4.74001741409302
0.8 4.74584674835205
0.82 4.69214582443237
0.84 4.70210266113281
0.86 4.69044017791748
0.88 4.71248483657837
0.9 4.766685962677
0.92 4.75997495651245
0.940000000000001 4.72668933868408
0.960000000000001 4.71324777603149
0.980000000000001 4.71860074996948
1 4.70157289505005
};
\addplot [semithick, black, dashed]
table {%
0 4.58966101118138
0.02 4.58752889141479
0.04 4.5864369641871
0.06 4.58637812293328
0.08 4.58980324561968
0.1 4.59073430595159
0.12 4.5794989874436
0.14 4.58275761739179
0.16 4.58070245378178
0.18 4.5655214154823
0.2 4.56156332474374
0.22 4.58109190112159
0.24 4.56027809246298
0.26 4.57142463241991
0.28 4.56628042982638
0.3 4.5841636768917
0.32 4.60015352220269
0.34 4.60539097746522
0.36 4.623795532681
0.38 4.63697801772565
0.4 4.62926810812352
0.42 4.62205931208351
0.44 4.59185886307324
0.46 4.60382516041881
0.48 4.63249140784437
0.5 4.65229545740424
0.52 4.66451264466782
0.54 4.67868608535656
0.56 4.71494178689265
0.58 4.74444150803121
0.6 4.73958008117334
0.62 4.72810550016589
0.64 4.78279924629234
0.66 4.80749265840588
0.68 4.8147370402964
0.7 4.79568127605334
0.72 4.79208937506806
0.74 4.84446063691436
0.76 4.81685398106912
0.78 4.76327837968624
0.8 4.76457957271389
0.82 4.70246698615537
0.84 4.71077552896357
0.86 4.6947133027551
0.88 4.71687700431492
0.9 4.7741259266026
0.92 4.76357110582256
0.940000000000001 4.72693643181042
0.960000000000001 4.71339805084577
0.980000000000001 4.71515653536922
1 4.69632364684386
};
\addplot [semithick, blue]
table {%
0 4.59025430679321
0.02 4.59063005447388
0.04 4.59184408187866
0.06 4.59139394760132
0.08 4.5915904045105
0.1 4.59266424179077
0.12 4.58828496932983
0.14 4.58995199203491
0.16 4.588050365448
0.18 4.5844669342041
0.2 4.5824933052063
0.22 4.5835337638855
0.24 4.58627080917358
0.26 4.58630132675171
0.28 4.58338356018066
0.3 4.58425378799438
0.32 4.58466482162476
0.34 4.60086250305176
0.36 4.60010385513306
0.38 4.60607576370239
0.4 4.63201427459717
0.42 4.63472127914429
0.44 4.64281320571899
0.46 4.65605592727661
0.48 4.64614868164062
0.5 4.62841033935547
0.52 4.67426109313965
0.54 4.66733598709106
0.56 4.67102003097534
0.58 4.66841793060303
0.6 4.67116212844849
0.62 4.64714479446411
0.64 4.63419818878174
0.66 4.61900520324707
0.68 4.64819478988647
0.7 4.64553260803223
0.72 4.65133380889893
0.74 4.63320589065552
0.76 4.63718318939209
0.78 4.62267017364502
0.8 4.64557504653931
0.82 4.61712598800659
0.84 4.64768981933594
0.86 4.64028882980347
0.88 4.64189386367798
0.9 4.67032527923584
0.92 4.63322734832764
0.940000000000001 4.60807847976685
0.960000000000001 4.63029050827026
0.980000000000001 4.61804437637329
1 4.60655450820923
};
\addplot [semithick, black, dashed]
table {%
0 4.58966101118138
0.02 4.58955956081935
0.04 4.59735467867489
0.06 4.59366074057664
0.08 4.59863981247258
0.1 4.60260863669828
0.12 4.58614785713916
0.14 4.59278348250191
0.16 4.58728012171564
0.18 4.57991100914683
0.2 4.58324800975921
0.22 4.5832894829115
0.24 4.59304929097869
0.26 4.59516633627534
0.28 4.58882389580392
0.3 4.59223669698128
0.32 4.592495399419
0.34 4.62462190074321
0.36 4.6195680061298
0.38 4.62657263293918
0.4 4.66906553731396
0.42 4.67377865586779
0.44 4.68722622374622
0.46 4.70567296062533
0.48 4.6839501213803
0.5 4.65334641144492
0.52 4.71926233177082
0.54 4.70801212643691
0.56 4.71221252246761
0.58 4.70715307408496
0.6 4.70990364266556
0.62 4.67881859317402
0.64 4.66123279709969
0.66 4.63956180020859
0.68 4.67566685958974
0.7 4.67283516407365
0.72 4.677473634144
0.74 4.65428232994138
0.76 4.65721402136129
0.78 4.64163079662163
0.8 4.66765875141699
0.82 4.63144991763689
0.84 4.66341070224677
0.86 4.65401989873739
0.88 4.65327445285068
0.9 4.68235362337752
0.92 4.64401497671768
0.940000000000001 4.61767283517204
0.960000000000001 4.64343723422209
0.980000000000001 4.62981876098172
1 4.61668154161746
};
\addplot [semithick, blue]
table {%
0 4.59025430679321
0.02 4.58968639373779
0.04 4.58959197998047
0.06 4.59257793426514
0.08 4.59150505065918
0.1 4.58605003356934
0.12 4.58758020401001
0.14 4.59014368057251
0.16 4.58614540100098
0.18 4.58794832229614
0.2 4.58283853530884
0.22 4.57947301864624
0.24 4.58197259902954
0.26 4.58751630783081
0.28 4.5897479057312
0.3 4.59646368026733
0.32 4.60616540908813
0.34 4.61001968383789
0.36 4.61512136459351
0.38 4.61183023452759
0.4 4.61269569396973
0.42 4.63634443283081
0.44 4.63463163375854
0.46 4.62416744232178
0.48 4.6163501739502
0.5 4.60310173034668
0.52 4.60844421386719
0.54 4.59814929962158
0.56 4.58304166793823
0.58 4.59414148330688
0.6 4.60056495666504
0.62 4.59687852859497
0.64 4.5732364654541
0.66 4.52879476547241
0.68 4.51804304122925
0.7 4.5476245880127
0.72 4.54501104354858
0.74 4.52476739883423
0.76 4.52753877639771
0.78 4.56591272354126
0.8 4.52180814743042
0.82 4.53791952133179
0.84 4.52125215530396
0.86 4.48923492431641
0.88 4.53060483932495
0.9 4.50458574295044
0.92 4.56839466094971
0.940000000000001 4.60252904891968
0.960000000000001 4.57530546188354
0.980000000000001 4.56523036956787
1 4.6044487953186
};
\addplot [semithick, black, dashed]
table {%
0 4.58966101118138
0.02 4.58773495610415
0.04 4.58910377950958
0.06 4.60268749836
0.08 4.60311273759488
0.1 4.58514102836428
0.12 4.59229955046153
0.14 4.59979832160644
0.16 4.58534133720425
0.18 4.59141772776058
0.2 4.57908769595284
0.22 4.5716021240034
0.24 4.57744941384649
0.26 4.59275179687954
0.28 4.59368880420751
0.3 4.60734802295849
0.32 4.62608934813349
0.34 4.62945883021159
0.36 4.63985473937878
0.38 4.63253472980438
0.4 4.63036248510412
0.42 4.6708136313589
0.44 4.66724248358048
0.46 4.64836444829033
0.48 4.63415067714578
0.5 4.61161636134827
0.52 4.61590164417068
0.54 4.60001345462341
0.56 4.57870396930181
0.58 4.59108987678945
0.6 4.59811028170646
0.62 4.59256704016279
0.64 4.56199863067005
0.66 4.50456629470547
0.68 4.49205453626681
0.7 4.52916496575902
0.72 4.52752208687736
0.74 4.50531499286219
0.76 4.51009716461574
0.78 4.55565471811837
0.8 4.5048624421895
0.82 4.52142048384382
0.84 4.50578446670394
0.86 4.47269943904501
0.88 4.52110453256042
0.9 4.49476796693203
0.92 4.56355762037902
0.940000000000001 4.5981780123073
0.960000000000001 4.57011696364627
0.980000000000001 4.55990991223033
1 4.59944658292746
};
\addplot [semithick, blue]
table {%
0 4.59025430679321
0.02 4.58993291854858
0.04 4.58828687667847
0.06 4.58673715591431
0.08 4.58647584915161
0.1 4.58679914474487
0.12 4.58093309402466
0.14 4.57858324050903
0.16 4.58004522323608
0.18 4.58069467544556
0.2 4.57732391357422
0.22 4.58307886123657
0.24 4.57779598236084
0.26 4.57541656494141
0.28 4.57657051086426
0.3 4.57395505905151
0.32 4.57032060623169
0.34 4.5658860206604
0.36 4.57563924789429
0.38 4.56890535354614
0.4 4.57318639755249
0.42 4.5772442817688
0.44 4.57462644577026
0.46 4.54176807403564
0.48 4.54912614822388
0.5 4.55529403686523
0.52 4.52230501174927
0.54 4.53235387802124
0.56 4.54734945297241
0.58 4.53350400924683
0.6 4.52119970321655
0.62 4.53648662567139
0.64 4.5282039642334
0.66 4.53623723983765
0.68 4.53000974655151
0.7 4.54029369354248
0.72 4.52941608428955
0.74 4.56820154190063
0.76 4.56011962890625
0.78 4.58147573471069
0.8 4.57825136184692
0.82 4.6007022857666
0.84 4.5929708480835
0.86 4.55044555664062
0.88 4.54803228378296
0.9 4.54284381866455
0.92 4.54278135299683
0.940000000000001 4.50603437423706
0.960000000000001 4.47492980957031
0.980000000000001 4.46057891845703
1 4.46184492111206
};
\addplot [semithick, black, dashed]
table {%
0 4.58966101118138
0.02 4.5913798622797
0.04 4.58547114401067
0.06 4.57936402276749
0.08 4.58042865365809
0.1 4.57960541256741
0.12 4.55919598616163
0.14 4.55461684452769
0.16 4.55522280022292
0.18 4.55638791467894
0.2 4.54763096871673
0.22 4.56569790671555
0.24 4.55514800376547
0.26 4.55062647900699
0.28 4.55179034864484
0.3 4.54844632105264
0.32 4.54338005193113
0.34 4.53614767560103
0.36 4.5585802572609
0.38 4.54941622287838
0.4 4.56126970252752
0.42 4.56745470105095
0.44 4.56554357900661
0.46 4.51413353813937
0.48 4.52767059448891
0.5 4.53577930267083
0.52 4.48541948438417
0.54 4.50283873526368
0.56 4.52490925730818
0.58 4.50793784678163
0.6 4.49062482305273
0.62 4.51686197629047
0.64 4.50822435723153
0.66 4.51852791992199
0.68 4.51139339557582
0.7 4.52421535356426
0.72 4.51311859727844
0.74 4.56048705612228
0.76 4.54812134336671
0.78 4.57492368142619
0.8 4.56965734064883
0.82 4.59648342541589
0.84 4.58880858923044
0.86 4.54160159656132
0.88 4.53871491108091
0.9 4.53367811843094
0.92 4.53056659142015
0.940000000000001 4.49321717751754
0.960000000000001 4.4622990541723
0.980000000000001 4.44941668560612
1 4.45319281817539
};

\nextgroupplot[
tick align=outside,
tick pos=left,
x grid style={darkgray176},
xmin=-0.05, xmax=1.05,
xtick style={color=black},
y grid style={darkgray176},
ymin=4.22494469146819, ymax=4.88066004894983,
ytick style={color=black}
]
\addplot [semithick, red]
table {%
0 4.5917444229126
0.02 4.5931282043457
0.04 4.59003734588623
0.06 4.58871746063232
0.08 4.58617305755615
0.1 4.59027576446533
0.12 4.58533382415771
0.14 4.58408355712891
0.16 4.58474588394165
0.18 4.58715152740479
0.2 4.58787488937378
0.22 4.59107971191406
0.24 4.58869934082031
0.26 4.58302927017212
0.28 4.58432817459106
0.3 4.59433174133301
0.32 4.58774471282959
0.34 4.59190559387207
0.36 4.60057067871094
0.38 4.57638549804688
0.4 4.57976770401001
0.42 4.57401847839355
0.44 4.58245277404785
0.46 4.6015362739563
0.48 4.59027576446533
0.5 4.56962966918945
0.52 4.57221460342407
0.54 4.56937885284424
0.56 4.55707359313965
0.58 4.53658103942871
0.6 4.53151416778564
0.62 4.51104211807251
0.64 4.5061206817627
0.66 4.4684853553772
0.68 4.4547233581543
0.7 4.48142051696777
0.72 4.47336673736572
0.74 4.4259295463562
0.76 4.3993763923645
0.78 4.39080142974854
0.8 4.38441467285156
0.82 4.40545511245728
0.84 4.39514684677124
0.86 4.38735866546631
0.88 4.38368225097656
0.9 4.38507556915283
0.92 4.38021087646484
0.940000000000001 4.32068920135498
0.960000000000001 4.32121515274048
0.980000000000001 4.27816581726074
1 4.26010513305664
};
\addplot [semithick, black, dashed]
table {%
0 4.58966101118138
0.02 4.59268717723046
0.04 4.58381363243605
0.06 4.58476147891541
0.08 4.5772737459526
0.1 4.59412229118622
0.12 4.5777661702421
0.14 4.57638524505657
0.16 4.57358566528419
0.18 4.58614978184785
0.2 4.59099243726389
0.22 4.60115939034652
0.24 4.59544210565938
0.26 4.58728805771992
0.28 4.5909381423137
0.3 4.61556403800043
0.32 4.59912845586912
0.34 4.60887807690182
0.36 4.62348909913715
0.38 4.57507450343685
0.4 4.5812093759618
0.42 4.56774166001864
0.44 4.58220488001283
0.46 4.61505559437336
0.48 4.59794740140478
0.5 4.56350905813169
0.52 4.56984471723566
0.54 4.5682463745631
0.56 4.54941174349783
0.58 4.51986892804467
0.6 4.51383492571129
0.62 4.48562979928835
0.64 4.4803137477497
0.66 4.43156961641438
0.68 4.41561028378318
0.7 4.45471137510068
0.72 4.44715979538322
0.74 4.38861239821717
0.76 4.35924111503422
0.78 4.35209018296526
0.8 4.34938676969589
0.82 4.37758983491855
0.84 4.37274802992552
0.86 4.36671691415653
0.88 4.36615063589296
0.9 4.36988546636657
0.92 4.36881881351852
0.940000000000001 4.30991951867014
0.960000000000001 4.31285998740165
0.980000000000001 4.27116402260826
1 4.25474993499008
};
\addplot [semithick, red]
table {%
0 4.5917444229126
0.02 4.59118366241455
0.04 4.59088706970215
0.06 4.59081840515137
0.08 4.59318733215332
0.1 4.59350776672363
0.12 4.58934688568115
0.14 4.59178781509399
0.16 4.59142875671387
0.18 4.5861439704895
0.2 4.58348274230957
0.22 4.58854484558105
0.24 4.58045196533203
0.26 4.58597040176392
0.28 4.58517980575562
0.3 4.59103488922119
0.32 4.59871387481689
0.34 4.60077571868896
0.36 4.61053037643433
0.38 4.61835813522339
0.4 4.61440658569336
0.42 4.61231899261475
0.44 4.59576320648193
0.46 4.60337543487549
0.48 4.62202167510986
0.5 4.63576126098633
0.52 4.64402675628662
0.54 4.65364646911621
0.56 4.68005847930908
0.58 4.70386838912964
0.6 4.70326900482178
0.62 4.69523620605469
0.64 4.73754596710205
0.66 4.75831604003906
0.68 4.76421928405762
0.7 4.75340175628662
0.72 4.75197982788086
0.74 4.79793453216553
0.76 4.77880668640137
0.78 4.73658990859985
0.8 4.74197816848755
0.82 4.69115447998047
0.84 4.69947242736816
0.86 4.69001960754395
0.88 4.71199750900269
0.9 4.76675796508789
0.92 4.76044273376465
0.940000000000001 4.72952461242676
0.960000000000001 4.71860933303833
0.980000000000001 4.72216987609863
1 4.70867156982422
};
\addplot [semithick, black, dashed]
table {%
0 4.58966101118138
0.02 4.58752883257261
0.04 4.58643997925971
0.06 4.5863882624511
0.08 4.58982751574787
0.1 4.59077643652918
0.12 4.57955854713733
0.14 4.58284382479854
0.16 4.58082265102315
0.18 4.56568027097167
0.2 4.56176534153053
0.22 4.581350957763
0.24 4.56059907704027
0.26 4.57181918598431
0.28 4.56676154867808
0.3 4.58473429167341
0.32 4.60082499482539
0.34 4.60617167326713
0.36 4.6247173473484
0.38 4.63804881278228
0.4 4.63052592733557
0.42 4.62352268616972
0.44 4.59355452335519
0.46 4.60575191478632
0.48 4.63465261796608
0.5 4.65476658901219
0.52 4.66725707851751
0.54 4.68171597173531
0.56 4.71822015936191
0.58 4.74800805382048
0.6 4.74348728562183
0.62 4.73236912425886
0.64 4.78732425854795
0.66 4.81237577176292
0.68 4.81994515386485
0.7 4.80140308413132
0.72 4.79822849863005
0.74 4.85085480542794
0.76 4.82367209875598
0.78 4.77061894226535
0.8 4.77244013372932
0.82 4.71101778163849
0.84 4.71980199386432
0.86 4.70432892457996
0.88 4.72694099149834
0.9 4.7844748152406
0.92 4.77453859516084
0.940000000000001 4.73847236682007
0.960000000000001 4.72549904985644
0.980000000000001 4.72795851656927
1 4.7098503517713
};
\addplot [semithick, red]
table {%
0 4.5917444229126
0.02 4.59224891662598
0.04 4.59415817260742
0.06 4.59210586547852
0.08 4.59249210357666
0.1 4.59305477142334
0.12 4.58934354782104
0.14 4.59231090545654
0.16 4.59101009368896
0.18 4.58767938613892
0.2 4.58782577514648
0.22 4.58849620819092
0.24 4.59196662902832
0.26 4.59107208251953
0.28 4.58885097503662
0.3 4.58933782577515
0.32 4.58974838256836
0.34 4.60701656341553
0.36 4.60604667663574
0.38 4.60955238342285
0.4 4.63425445556641
0.42 4.63914728164673
0.44 4.64964771270752
0.46 4.66495418548584
0.48 4.65212965011597
0.5 4.63380432128906
0.52 4.6779088973999
0.54 4.67138385772705
0.56 4.67575979232788
0.58 4.67324304580688
0.6 4.6782283782959
0.62 4.65547752380371
0.64 4.64504432678223
0.66 4.63040828704834
0.68 4.65918064117432
0.7 4.65677213668823
0.72 4.66565942764282
0.74 4.64374446868896
0.76 4.64719104766846
0.78 4.63444900512695
0.8 4.65857744216919
0.82 4.62845945358276
0.84 4.66037559509277
0.86 4.65496349334717
0.88 4.65690422058105
0.9 4.68607091903687
0.92 4.65119028091431
0.940000000000001 4.62872982025146
0.960000000000001 4.6527738571167
0.980000000000001 4.6430492401123
1 4.62937545776367
};
\addplot [semithick, black, dashed]
table {%
0 4.58966101118138
0.02 4.58955968547814
0.04 4.59735846811758
0.06 4.59367331669433
0.08 4.598663255348
0.1 4.60265265916297
0.12 4.58621416214009
0.14 4.59288024732227
0.16 4.58741315660454
0.18 4.58008692486529
0.2 4.58346570808413
0.22 4.58356082467906
0.24 4.59338457344752
0.26 4.59558011769498
0.28 4.58931339892501
0.3 4.5928160877063
0.32 4.59319296106983
0.34 4.6254587533019
0.36 4.62055803455986
0.38 4.62773924445047
0.4 4.67041191197766
0.42 4.67533205887077
0.44 4.68897383705537
0.46 4.70766139801646
0.48 4.68619958420982
0.5 4.65589506003753
0.52 4.72204843592294
0.54 4.71108966151503
0.56 4.71560164728779
0.58 4.71083148045038
0.6 4.71394758498418
0.62 4.68322641723504
0.64 4.66599868166428
0.66 4.64477234339357
0.68 4.68125395167688
0.7 4.67877912137685
0.72 4.68375970672607
0.74 4.66112706063751
0.76 4.66445836097296
0.78 4.64936778943408
0.8 4.67588809005936
0.82 4.64029082417141
0.84 4.67275689385242
0.86 4.66392633754137
0.88 4.66381674483051
0.9 4.6933229369344
0.92 4.65565562829565
0.940000000000001 4.63015273148038
0.960000000000001 4.65635533547474
0.980000000000001 4.64349285222437
1 4.63116421068632
};
\addplot [semithick, red]
table {%
0 4.5917444229126
0.02 4.5919623374939
0.04 4.59256362915039
0.06 4.59590768814087
0.08 4.59580230712891
0.1 4.59077262878418
0.12 4.59282445907593
0.14 4.59590625762939
0.16 4.59159278869629
0.18 4.59321594238281
0.2 4.58861303329468
0.22 4.58576583862305
0.24 4.58868408203125
0.26 4.59629917144775
0.28 4.59826517105103
0.3 4.60501527786255
0.32 4.61378192901611
0.34 4.61583662033081
0.36 4.6195592880249
0.38 4.61635684967041
0.4 4.61664772033691
0.42 4.640793800354
0.44 4.6397876739502
0.46 4.62862348556519
0.48 4.61957740783691
0.5 4.60538768768311
0.52 4.6094799041748
0.54 4.59878349304199
0.56 4.58400821685791
0.58 4.59548807144165
0.6 4.60193157196045
0.62 4.59750366210938
0.64 4.57487392425537
0.66 4.53055000305176
0.68 4.51900386810303
0.7 4.54843235015869
0.72 4.54667806625366
0.74 4.52529382705688
0.76 4.52745342254639
0.78 4.56721305847168
0.8 4.52198839187622
0.82 4.53518486022949
0.84 4.52029895782471
0.86 4.49094438552856
0.88 4.53347301483154
0.9 4.50759553909302
0.92 4.57192039489746
0.940000000000001 4.60410594940186
0.960000000000001 4.57838535308838
0.980000000000001 4.5695276260376
1 4.60812473297119
};
\addplot [semithick, black, dashed]
table {%
0 4.58966101118138
0.02 4.58773496359456
0.04 4.58910648630282
0.06 4.60269853339841
0.08 4.60313657170637
0.1 4.58517920627964
0.12 4.59235864475574
0.14 4.59988354179181
0.16 4.58546059273696
0.18 4.59157384195547
0.2 4.57928753621347
0.22 4.57185564803354
0.24 4.5777690599519
0.26 4.59316018476077
0.28 4.59417988441598
0.3 4.60793237489115
0.32 4.62680433100008
0.34 4.63029947352593
0.36 4.64085077758456
0.38 4.6336946415956
0.4 4.63170563903837
0.42 4.67236359353769
0.44 4.66897524055949
0.46 4.65030994197542
0.48 4.63633173659046
0.5 4.6140237346032
0.52 4.61854685944705
0.54 4.60291699297263
0.56 4.58197833921066
0.58 4.59466771008879
0.6 4.60200779032121
0.62 4.59684479999503
0.64 4.56670644450978
0.66 4.50972749836846
0.68 4.49777895882681
0.7 4.53528509730703
0.72 4.53411601313238
0.74 4.51238674812116
0.76 4.51761487440527
0.78 4.56358161557418
0.8 4.51361996202131
0.82 4.53068452244444
0.84 4.51558686231157
0.86 4.48308472525178
0.88 4.53196400442213
0.9 4.50655781476387
0.92 4.57563368123585
0.940000000000001 4.6109864225499
0.960000000000001 4.58395466842958
0.980000000000001 4.57465453240115
1 4.6146975181583
};
\addplot [semithick, red]
table {%
0 4.5917444229126
0.02 4.5915060043335
0.04 4.59047508239746
0.06 4.5896258354187
0.08 4.58951091766357
0.1 4.58953189849854
0.12 4.583083152771
0.14 4.58047819137573
0.16 4.58087205886841
0.18 4.57951736450195
0.2 4.57584095001221
0.22 4.58022260665894
0.24 4.57547664642334
0.26 4.57248592376709
0.28 4.57347679138184
0.3 4.57118463516235
0.32 4.56657314300537
0.34 4.56325626373291
0.36 4.57311248779297
0.38 4.56789302825928
0.4 4.5731258392334
0.42 4.57607746124268
0.44 4.57394313812256
0.46 4.54153680801392
0.48 4.54674291610718
0.5 4.55205965042114
0.52 4.51896953582764
0.54 4.52853775024414
0.56 4.5443754196167
0.58 4.53185081481934
0.6 4.52055501937866
0.62 4.53540229797363
0.64 4.52727270126343
0.66 4.53583574295044
0.68 4.53087282180786
0.7 4.54232215881348
0.72 4.53322458267212
0.74 4.57140254974365
0.76 4.56161165237427
0.78 4.58506155014038
0.8 4.57932806015015
0.82 4.60266971588135
0.84 4.59658622741699
0.86 4.55602598190308
0.88 4.55487442016602
0.9 4.5494532585144
0.92 4.54636383056641
0.940000000000001 4.50939083099365
0.960000000000001 4.48014736175537
0.980000000000001 4.46651935577393
1 4.47117805480957
};
\addplot [semithick, black, dashed]
table {%
0 4.58966101118138
0.02 4.5913795770118
0.04 4.58547410137321
0.06 4.57937313683296
0.08 4.58044666383273
0.1 4.57963967440258
0.12 4.55924925328224
0.14 4.55469042056494
0.16 4.55532758376173
0.18 4.55652737042165
0.2 4.54781128266762
0.22 4.56593163080325
0.24 4.55543935768629
0.26 4.55098385248005
0.28 4.55223922076897
0.3 4.54898699542731
0.32 4.54402708947324
0.34 4.53689132482004
0.36 4.55947010132069
0.38 4.55043018439847
0.4 4.56244204007614
0.42 4.56880418809778
0.44 4.56709863151056
0.46 4.51588676844945
0.48 4.52966053537577
0.5 4.53802713327476
0.52 4.48788025160011
0.54 4.50558086853813
0.56 4.5279679171837
0.58 4.51131551575166
0.6 4.49438902086361
0.62 4.52096167571498
0.64 4.51268228221872
0.66 4.52339300808452
0.68 4.51666822637802
0.7 4.5298407208378
0.72 4.51926116402849
0.74 4.5670158460774
0.76 4.55517534125505
0.78 4.58245451711115
0.8 4.57763902534721
0.82 4.60498921003849
0.84 4.59792750306555
0.86 4.55143740580866
0.88 4.54934551112062
0.9 4.54503260401405
0.92 4.5424966541079
0.940000000000001 4.50620167253508
0.960000000000001 4.47644537844205
0.980000000000001 4.46459410571718
1 4.46947525855592
};
\end{groupplot}
\node[anchor=north] (title-x) at ($(group c1r2.south east)!0.5!(group c2r2.south west)-(0,0.5cm)$) {\Large $t$};
\node[anchor=south, rotate=90] (title-y) at ($(group c1r1.south west)!0.5!(group c1r2.north west)-(0.85,0cm)$) {\Large $ Y_t = u(t,X_t)$};
\node[anchor=south] (title-z1) at ($(group c1r1.north east)!0.5!(group c2r1.north west)-(0.6,0cm)$) {\Large Using 20k iterations};
\node[anchor=south] (title-z2) at ($(group c1r2.north east)!0.5!(group c2r2.north west)-(0.6,0cm)$) {\Large Using 50k iterations};
\path ($(group c1r1.south east) -(0,0.7cm) $) -- node[below]{\ref{zelda}} ($(group c2r1.south west)-(0,0.7cm)$);
\end{tikzpicture}

%% file: figures/100D_problem/shifted_example.tex
% This file was created with tikzplotlib v0.10.1.
\begin{tikzpicture}

\definecolor{darkgray176}{RGB}{176,176,176}
\definecolor{darkorange25512714}{RGB}{255,127,14}
\definecolor{forestgreen4416044}{RGB}{44,160,44}
\definecolor{steelblue31119180}{RGB}{31,119,180}

\begin{axis}[
tick align=outside,
tick pos=left,
title={generated trajectories},
x grid style={darkgray176},
xlabel={dim 1},
xmin=-0.5, xmax=3.5,
xtick style={color=black},
y grid style={darkgray176},
ylabel={dim 2},
ymin=-0.5, ymax=3.5,
ytick style={color=black}
]
\addplot [very thick, darkorange25512714, mark=*, mark size=3, mark options={solid}]
table {%
0 0
0.217292606830597 0.130413621664047
0.496757715940475 0.705424427986145
0.677012801170349 0.752269268035889
0.83798748254776 0.645680487155914
0.990877687931061 1.09664440155029
0.908147692680359 1.10867810249329
0.945754766464233 1.36814391613007
1.29756140708923 1.23969268798828
1.57836067676544 1.61896204948425
1.59543883800507 1.54563045501709
1.76633250713348 1.93290197849274
2.15237045288086 2.16918420791626
2.10399794578552 2.31360220909119
2.24111151695251 2.25242829322815
2.3073992729187 2.50631332397461
2.73518300056458 2.7643358707428
2.59736847877502 2.99291563034058
2.65514302253723 3.11496806144714
2.81012439727783 3.24492216110229
3.09260272979736 3.22308373451233
};
\addplot [very thick, steelblue31119180, mark=*, mark size=3, mark options={solid}]
table {%
0 0
-0.0629750639200211 0.189047455787659
0.145277470350266 0.288629800081253
0.256100475788116 0.482483506202698
0.140017986297607 0.22880083322525
0.10454323142767 0.624352693557739
0.0855461433529854 0.501739621162415
0.149085268378258 0.615288555622101
0.240163117647171 0.767765700817108
0.474654376506805 0.876863777637482
0.435459077358246 1.01756608486176
0.41266068816185 0.904009222984314
0.503550052642822 1.29607176780701
0.394577741622925 1.18336272239685
0.501839756965637 1.22336184978485
0.803385615348816 1.15171980857849
0.761094748973846 1.23733389377594
0.88884449005127 1.47220051288605
0.878366470336914 1.49727618694305
0.937044262886047 1.323317527771
0.988030552864075 1.67754805088043
};
\addplot [semithick, forestgreen4416044, mark=x, mark size=5, mark options={solid}, only marks]
table {%
3 3
};
\end{axis}

\end{tikzpicture}